\newif\ifprint
\printfalse
\documentclass[twoside,10pt]{NicePaperCMSS}

\usepackage{amsmath}
\usepackage{amsthm}
\usepackage{amssymb}
\usepackage{amsfonts}
\usepackage{amscd}
\usepackage{wasysym}
\usepackage{latexsym}
\usepackage{bm}
\usepackage{hyperref}
\usepackage[utf8]{inputenc}
\usepackage{textcomp}
\usepackage{color}
\usepackage[pdftex]{graphicx}
\usepackage{titlesec}
\usepackage{xspace}
\usepackage{pdfsync}
\usepackage[numbers,square,comma]{natbib}
\usepackage{fixltx2e} %for \textsubscript

\ifprint
	\definecolor{linkred}{rgb}{0,0,0} % black
	\definecolor{linkblue}{rgb}{0,0,0} % black
\else
	\definecolor{linkred}{rgb}{0.7,0.2,0.2}
	\definecolor{linkblue}{rgb}{0,0.2,0.6}
\fi

\numberwithin{equation}{section} % standardize numbering

\catcode`@=11 \baselineskip 4.5mm
\def\ps@handbook{\def\@oddhead{\hfill \leftmark \hfill\thepage }
\def\@evenhead{\thepage \hfill \rightmark \hfill}
\def\@oddfoot{}
\def\@evenfoot{}}
\def\@evenhead{}
\def\@oddfoot{}
\def\@evenfoot{\hfill}
\def\list#1#2{\ifnum \@listdepth >5\relax \@toodeep \else \global
\advance \@listdepth\@ne \fi \rightmargin \z@ \listparindent\z@
\itemindent\z@ \csname @list\romannumeral\the\@listdepth\endcsname
\def\@itemlabel{#1}\let\makelabel\@mklab \@nmbrlistfalse #2\relax
\@trivlist \parskip -\parsep \parindent\listparindent \advance
\linewidth -\rightmargin \advance\linewidth -\leftmargin \advance
\@totalleftmargin \leftmargin \parshape \@ne \@totalleftmargin
\linewidth \ignorespaces}
\catcode`@=12
\renewcommand\thesection{\arabic{section}}
\renewcommand\thesubsection{\arabic{subsection}}
\renewcommand\thesubsubsection{\arabic{subsubsection}}
\renewcommand{\theequation}{\thesection.\arabic{equation}}

\def\thebibliography#1{\section*{References}
\list{[\arabic{enumi}]}{\settowidth \labelwidth{[#1]} \leftmargin
\labelwidth \advance \leftmargin \labelsep \usecounter{enumi}}
\def\newblock{\hskip .11em plus .33em minus .07em} \sloppy
\clubpenalty 4000 \widowpenalty 4000 \sfcode`\.=1000 \relax}

\titleformat{\section}{\normalfont\large\bfseries}{\thesection.}{0.5em}{}
\titleformat{\subsection}{\normalfont\bfseries}{\thesection.\thesubsection.}{0.5em}{}
\titleformat{\subsubsection}[runin]{\normalfont\bfseries}{\thesection.\thesubsection.\thesubsubsection.}{0.5em}{}[\kern0.5em]%IM

\newtheorem{theorem}[equation]{Theorem}
\newtheorem{proposition}[equation]{Proposition}
\newtheorem{corollary}[equation]{Corollary}

\newtheorem{lemma}[equation]{Lemma}
\newtheorem{definition}[equation]{Definition}
\theoremstyle{remark}
\newtheorem{remark}[equation]{Remark}

\newcommand{\bdoc}{
\begin{document}}
\newcommand{\edoc}{\end{document}}

\newcommand{\bcent}{\begin{center}}
\newcommand{\ecent}{\end{center}}

\newcommand{\benum}{\begin{enumerate}}
\newcommand{\eenum}{\end{enumerate}}
\newcommand{\bitem}{\begin{itemize}}
\newcommand{\eitem}{\end{itemize}}

\newcommand{\btab}{\begin{tabular}}
\newcommand{\etab}{\end{tabular}}
\newcommand{\beqn}{\begin{eqnarray}}
\newcommand{\eeqn}{\end{eqnarray}}

\newcommand{\bmath}{\begin{math}}
\newcommand{\emath}{\end{math}}

\newcommand{\noin}{\noindent}

\providecommand{\tb}[1]{\textbf{#1}}

\newcommand{\bsh}{\backslash}

\newcommand{\ds}{\displaystyle}

\newcommand{\ub}{\underbrace}
\newcommand{\ob}{\overbrace}

\providecommand{\F}[1]{\mathbb{#1}}  %change back to \mathbb for black board number fields

\newcommand{\FF}{\F F}
\providecommand{\Fn}[1]{\FF_{#1}}
\newcommand{\Fp}{\Fn{p}}
\newcommand{\Fq}{\Fn{q}}
\newcommand{\Fpm}{\Fn{p^m}}
\newcommand{\Fpn}{\Fn{p^n}}
\newcommand{\Fpr}{\Fn{p^r}}

\newcommand{\ZZ}{\mathbf Z} %always black board even if \F is not
\providecommand{\Zn}[1]{\ZZ_{#1}}
\providecommand{\ZnZ}[1]{\ZZ/#1\ZZ}
\newcommand{\Zp}{\Zn{p}}
\newcommand{\ZpZ}{\ZnZ{p}}

\newcommand{\NN}{\F{N}}
\newcommand{\QQ}{\F{Q}}
\newcommand{\RR}{\F{R}}
\newcommand{\CC}{\mathbf C}
\newcommand{\QQbar}{\overline{\QQ}}
\newcommand{\Zero}{\mathbb{00}}

\newcommand{\EE}{\F E}
\newcommand{\II}{\F I}
\newcommand{\KK}{\F K}
\newcommand{\MM}{\F M}
\newcommand{\XX}{\F X}
\newcommand{\PP}{\mathbf P}
\newcommand{\FA}{\F A}
\newcommand{\LL}{\F L}
\newcommand{\HH}{\mathbb H} % changed to \mathbb from \F to force blackboard font for hypercohomology in all cases
\newcommand{\FS}{\F S}
\newcommand{\TT}{\F T} %parametrising scheme
\newcommand{\FSig}{\F\Sig}
\newcommand{\FDel}{\F\Del}

\providecommand{\E}[1]{\hat{\F{#1}}}
\newcommand{\EC}{\E{C}}

\newcommand{\bQ}{\mathbf{Q}}
\newcommand{\bP}{\mathbf{P}}

\newcommand{\ind}{\mbox{ind}}

\newcommand{\fx}{f(x)}
\newcommand{\gx}{g(x)}

\newcommand{\x}{^\star} %units
\newcommand{\xs}{^{~\star}} %units with space
\providecommand{\U}[1]{\left(#1\right)\x} %units with parentheses
\newcommand{\xt}{^\times} %units

\newcommand{\iso}{\simeq}
\newcommand{\plus}{\oplus}
\newcommand{\Plus}{\bigoplus}
\newcommand{\tensor}{\otimes}
\newcommand{\Tensor}{\bigotimes}
\newcommand{\inject}{\hookrightarrow}
\newcommand{\linject}{\hookleftarrow}
\newcommand{\surject}{\twoheadrightarrow}
\newcommand{\tri}{\vartriangleleft}

\renewcommand{\ker}{\mbox{ker}\;}
\newcommand{\Imf}{\mbox{im}~\;}
\newcommand{\img}{\mbox{im}\;}
\newcommand{\Hom}{\mbox{Hom}}
\newcommand{\Sym}{\mbox{Sym}}
\newcommand{\End}{\mbox{End}\,}
\newcommand{\Endz}{\mbox{End}_0}
\newcommand{\Id}{\mbox{Id}}
\newcommand{\rk}{\mbox{rk}\,}
\newcommand{\Pic}{\mbox{Pic}}
\newcommand{\Jac}{\mbox{Jac}}
\newcommand{\ch}{\mbox{ch}\,}
\newcommand{\td}{\mbox{td}\,}

\providecommand{\Gal}[1]{\mbox{Gal}(#1)}
\providecommand{\GAL}[2]{\mbox{Gal}(#1/#2)}
\providecommand{\Sub}[1]{\mbox{Sub}(#1)}
\providecommand{\Lat}[1]{\mbox{Lat}(#1)}

\newcommand{\dup}{d_\wedge}
\newcommand{\drt}{d_>}

\newcommand{\MF}{\mathfrak}

\newcommand{\Div}{\,\,|\,\,}
\newcommand{\lcm}{\mbox{lcm}}

\providecommand{\leg}[2]{\left(\frac{#1}{#2}\right)}
\providecommand{\jac}[2]{\leg{#1}{#2}}
\providecommand{\qdc}[2]{\left[\frac{#1}{#2}\right]}

\newcommand{\w}{\omega}
\newcommand{\W}{\Omega}
\providecommand{\Cal}[1]{\mathcal{#1}}
\newcommand{\CL}{\Cal L}
\newcommand{\CO}{\Cal O}
\renewcommand{\O}{\Cal O}
\renewcommand{\o}{\Cal O}
\newcommand{\co}{\Cal O}
\newcommand{\ca}{\Cal A}
\newcommand{\CE}{\Cal E}
\newcommand{\CF}{\Cal F}
\newcommand{\CQ}{\Cal Q}
\newcommand{\CK}{\Cal K}
\newcommand{\CM}{\Cal M}
\newcommand{\CN}{\Cal N}
\newcommand{\CDee}{\Cal D}
\newcommand{\CP}{\Cal P}
\newcommand{\CS}{\Cal S}
\newcommand{\ClC}{\Cal C}
\newcommand{\CJ}{\Cal{J}}
\newcommand{\CA}{\Cal{A}}
\newcommand{\CB}{\Cal{B}}

\newcommand{\Nilp}{\mathbf{Nilp}}

\newcommand{\Si}{\Sigma}

\providecommand{\Ok}[1]{\CO_{#1}}
\newcommand{\OK}{\Ok{K}}

\renewcommand{\epsilon}{\varepsilon}
\newcommand{\ep}{\varepsilon}

\providecommand{\abs}[1]{\left|#1\right|}
\providecommand{\norm}[1]{\lVert#1\rVert}

\newcommand{\di}{\partial}
\providecommand{\ddy}[1]{\ds\frac{d}{d #1}}
\providecommand{\didiy}[1]{\ds\frac{\di}{\di #1}}
\newcommand{\ddx}{\ddy{x}}
\newcommand{\didix}{\didiy{x}}
\providecommand{\v}[1]{\vec{#1}}

\newcommand{\ihat}{\hat{\infty}}

\providecommand{\set}[1]{\left\{#1\right\}}
\providecommand{\lst}[1]{\ds\left[\,\,#1\,\,\right]}

\providecommand{\ip}[2]{\left( #1,#2\right)}
\providecommand{\IP}[2]{\left\langle #1,#2\right\rangle}
\providecommand{\bra}[1]{\left\langle\left. #1\right|\right.}
\providecommand{\ket}[1]{\left.\left| #1\right.\right\rangle}
\providecommand{\bkv}[4]{\left\langle\left.\tb{#1} #2\right|\tb{#3} #4\right\rangle}
\providecommand{\bk}[2]{\bkv{}{#1}{}{#2}}
\providecommand{\bkvop}[5]{\left\langle\tb{#1} #2\left| #5\right|\tb{#3} #4\right\rangle}
\providecommand{\bkop}[3]{\bkvop{}{#1}{}{#2}{#3}}
\providecommand{\comm}[2]{\left[#1,#2\right]}
\providecommand{\expn}[1]{\left\langle #1\right\rangle}
\providecommand{\gen}[1]{\expn{#1}}
\newcommand{\Del}{\Delta}
\newcommand{\Nab}{\nabla}
\newcommand{\Sig}{\Sigma}
\newcommand{\oline}{\overline}

\newcommand{\p}{\rho}
\providecommand{\sprod}[2]{\left\langle #1,#2\right\rangle}

\newcommand{\Ehat}{\overline E}
\newcommand{\Phihat}{\overline \Phi}
\newcommand{\phihat}{\overline \phi}

\newcommand{\QED}{\begin{flushright}\rule{2.5mm}{2.5mm}\\\end{flushright}}

\renewcommand{\gcd}{\mbox{gcd}}
\renewcommand{\deg}{\mbox{deg}}
\providecommand{\mb}[1]{\mathbf{#1}}
\newcommand{\h}{\mb h}
\newcommand{\V}{{V^\p_k}}

\definecolor{Blanc}{rgb}{1,1,1}
\definecolor{Violet}{rgb}{0.1,0,1}

%%%
%%%
%%%
%%%
%%%%%%%%%%%%%%%%%%%%%%%%%%%%%%%%%%%%%%%%%%%%%%%%%%%%%%%%%%%%%%%%%%%%%%%%%%%%%%%%%%%%%%%%%%%%%%%%%%%%
%%%%%%%%%%%%%%%%%%%%%%%%%%%%%%%%%%%%%%%%%%%%%%%%%%%%%%%%%%%%%%%%%%%%%%%%%%%%%%%%%%%%%%%%%%%%%%%%%%%%
%%%%%%%%%%%%%%%%%%%%%%%%%%%%%%%%%%%%%%%%%%%%%%%%%%%%%%%%%%%%%%%%%%%%%%%%%%%%%%%%%%%%%%%%%%%%%%%%%%%%

\setcounter{page}{1}
\long\def\replace#1{#1}

%    
%    Title
%
{\title[\fontfamily{cmss}\selectfont CONSTRUCTING CO-HIGGS BUNDLES ON \tb{CP}\textsuperscript2]{\fontfamily{cmss}\selectfont CONSTRUCTING CO-HIGGS BUNDLES ON \tb{CP}\textsuperscript2}    
\author{\fontfamily{cmss}\selectfont STEVEN RAYAN}
\address{\fontfamily{cmss}\selectfont Department of Mathematics, University of Toronto, 40 St. George Street, Toronto, Ontario, M5S 2E4~ CANADA}
\email{{\fontfamily{cmss}\selectfont rayan@math.toronto.edu}}

%    
%    Classification and abstract
%    

\begin{document}

\begin{abstract}
	
{On a complex manifold, a co-Higgs bundle is a holomorphic vector bundle with an endomorphism twisted by the tangent bundle.  The notion of generalized holomorphic bundle in Hitchin's generalized geometry coincides with that of co-Higgs bundle when the generalized complex manifold is ordinary complex.  Schwarzenberger's rank-2 vector bundle on the projective plane, constructed from a line bundle on the double cover $\CC\PP^1\times\CC\PP^1\to\CC\PP^2$, is naturally a co-Higgs bundle, with the twisted endomorphism, or ``Higgs field'', also descending from the double cover.  Allowing the branch conic to vary, we find that Schwarzenberger bundles give rise to an 8-dimensional moduli space of co-Higgs bundles.  After studying the deformation theory for co-Higgs bundles on complex manifolds, we conclude that a co-Higgs bundle arising from a Schwarzenberger bundle with nonzero Higgs field is rigid, in the sense that a nearby deformation is again Schwarzenberger.}

\end{abstract}  

\maketitle
\thispagestyle{empty}% suppress page number on title pages

\vspace{20pt}

\section{\fontfamily{cmss}\selectfont INTRODUCTION}

\noin  The goal of this paper is to explore an observation of Gualtieri \citep[][$\S$4.1]{MG:07}, namely that by considering ordinary complex manifolds in the context of generalized geometry, there is an enlargement of the category of holomorphic bundles.  The additional objects have the following form:\begin{definition}\label{DefnCoHiggs}  If $X$ is a complex manifold with tangent bundle $T_X$, then a \emph{\tb{co-Higgs bundle}} on $X$ is a holomorphic vector bundle $V\rightarrow X$ together with a map $\Phi\in H^0(X,(\emph{\End}V)\tensor T_X)$ for which $\Phi\wedge\Phi=0\in H^0(X,(\emph{\End V})\tensor\wedge^2T_X)$.\end{definition}\noin For generalized complex manifolds, there is an appropriate notion of bundle, called a ``generalized holomorphic bundle''.  When we consider ordinary complex manifolds as examples of generalized complex manifolds, the definition of generalized holomorphic bundle in \citep[][$\S$3.2]{MG:11} coincides with Definition \eqref{DefnCoHiggs}.  The key point is that a generalized holomorphic bundle on an ordinary complex manifold is not simply a holomorphic vector bundle, although holomorphic vector bundles are examples, arising when $\Phi=0$.

The reasoning behind the name ``co-Higgs bundle'' is that the object in Definition \eqref{DefnCoHiggs} resembles what is usually called a Higgs bundle, except that $\Phi$ takes values in $T_X^\vee$ for a Higgs bundle.  It is common in the theory of Higgs bundles to refer to the data $\Phi$ as a \emph{Higgs field} and we will use this terminology here.  We sometimes write \emph{integrable Higgs field} to emphasize that $\Phi$ satisfies $\Phi\wedge\Phi=0$.

Owing to the analogy with Higgs bundles, many features of Higgs bundles and their moduli spaces carry over to co-Higgs bundles.  For one, Higgs bundles come with a natural stability condition, discovered by Hitchin in \cite{NJH:86}, generalizing Mumford's slope stability for vector bundles.  We can adapt the stability condition for use here, allowing us to discuss moduli spaces of co-Higgs bundles.\begin{definition}A co-Higgs bundle $(V,\Phi)$ on a complex projective manifold $X$ is \emph{\tb{semistable}} if\beqn\frac{\emph{\deg}\,U}{\emph{\rk}U} & \leq & \frac{\emph{\deg}\,V}{\emph{\rk}V}\label{IneqStab}\eeqn\noin for all coherent subsheaves $0\neq U\subsetneq V$ satisfying $\Phi(U)\subseteq U\tensor T$, and \emph{\tb{stable}} if \eqref{IneqStab} is strict for all such $U$.\end{definition}\noin The projectivity assumption is used in the definition of the degree.  Stability for ordinary vector bundles without Higgs fields is recovered by taking $\Phi=0$.  When $V$ is fixed, we refer to $\Phi$ as (semi)stable whenever the pair $(V,\Phi)$ is (semi)stable.  There are situations in this paper where it will be necessary to consider pairs $(V,\phi)$ in which $\phi$ is an endomorphism taking values in a line bundle $L$.  For these objects, the stability condition is identical, simply with $L$ in place of $T_X$.

Examples of co-Higgs bundles have appeared in other studies recently.  In \cite{SSR:10}, we prove an existence theorem characterizing exactly those splitting types on $\CC\PP^1$ for which stable co-Higgs bundles exist.  There are no stable co-Higgs bundles with nonzero Higgs field on curves of genus $g>1$.  (When $g=1$, a co-Higgs bundle is the same thing as a Higgs bundle in the usual sense.)  In \cite{NJH:10IIa}, Hitchin constructs examples of generalized holomorphic bundles on complex manifolds.  A seed for the primary construction in our paper, using Schwarzenberger bundles, was planted in Hitchin's paper.

In the curve case, the only concern is stability.  Complex surfaces are a different story, at least at rank 2.  The integrability condition $\Phi\wedge\Phi=0$ is nontrivial, and is the main obstacle to finding examples of co-Higgs bundles.  The majority of this paper is occupied with constructing rank-2 examples on the complex projective plane.

For convenience, we will restrict from now on to trace-zero Higgs fields, which we signify by $\Phi\in H^0(X,(\Endz V)\tensor T_X)$.

\subsection{\fontfamily{cmss}\selectfont Results}

A portion of this paper is devoted to deformation theory.  As the most general part of the paper, it applies to co-Higgs bundles of arbitrary rank on complex manifolds of any dimension.  The condition $\Phi\wedge\Phi=0$ gives rise to the complex\beqn\Endz V\stackrel{\wedge\Phi}{\longrightarrow}(\Endz V)\tensor T_X\stackrel{\wedge\Phi}{\longrightarrow}(\Endz V)\tensor\wedge^2 T_X\stackrel{\wedge\Phi}{\longrightarrow}\cdots\nonumber\eeqn\noin that controls the deformation theory and whose spectral sequence computes the tangent space to the moduli space near $(V,\Phi)$.

Before constructing explicit examples, we justify the attention given to $\CC\PP^2$ amongst complex surfaces.  We provide a vanishing theorem (Theorem \ref{ThmVanishing}) for stable rank-2 co-Higgs bundles on general-type projective surfaces.  The proof of the theorem hints that stable rank-2 examples might be found wherever holomorphic sections of $\mb{S}^2T_X$ are plentiful.  The projective plane is a reasonable place to start, considering that $h^0(\mb{S}^2T_{\CC\PP^2})=27$.

If $\p$ is an irreducible element of $H^0(\CC\PP^2,\CO(2))$, then $\p=0$ defines a nonsingular conic, as well as a degree-2 covering of $\CC\PP^2$ by a smooth quadric $\CC\PP^1\times\CC\PP^1$ branched over the conic.  Each such conic also determines a sequence of rank-2 vector bundles, $\set{V^\p_k}_{k\ge0}$.  If $f^\p$ is the covering map, then we define\beqn V^\p_k:=f^\p_*\CO(0,k)\nonumber\eeqn\noin for each $k\geq0$, where $\CO(0,k):=\CO_1\tensor\CO_2(k)$, with $\CO_1$ pulled back from one ruling and $\CO_2$ from the other.  The bundles $V^\p_k$, called \emph{Schwarzenberger bundles}, were first studied in \cite{RLES:61I}.  They are likely the earliest examples of indecomposable holomorphic bundles on a complex surface. (Only $V^\p_0$ and $V^\p_1$ are decomposable.  For $k\geq2$, the bundles are indecomposable and stable.)  In Schwarzenberger's original study, the branch conic is a fixed nonsingular conic.  We allow the conic to vary and study the moduli problem for co-Higgs bundles whose underlying vector bundles are Schwarzenberger bundles arising from nonsingular conics.  

We show that Schwarzenberger bundles come naturally with integrable Higgs fields taking values in the tangent bundle of $\CC\PP^2$.  Like the Schwarzenberger bundles themselves, these Higgs fields descend from the line bundle on the double cover.  We show that for all $k$ the moduli space of co-Higgs bundles arising from Schwarzenberger bundles for nonsingular conics is 8-dimensional.  For $k=0$ and $k=2$, the moduli space is the total space of a vector bundle over a complex projective space.  For $k\geq3$, the moduli space admits two canonical fibrations, over projective spaces of different dimension.  For $k=1$, we construct a dense open set of the moduli space, which itself is the total space of a vector bundle.  For $k=0,1,2$, our moduli descriptions include the singular conics.

That the moduli space is 8-dimensional in every case is an application of the deformation theory, which we compute using cohomologies of exact sequences on the double cover.

\subsection{\fontfamily{cmss}\selectfont Range of the construction}

Let $H$ be the Chern class of $\CO_{\CC\PP^2}(1)$.  To measure the coverage our construction provides with respect to possible Chern classes $(c_1,c_2)\in\ZZ H\times\ZZ H^2$, we normalize the classes so that $(c_1,c_2)\in\set{0,-H}\times\ZZ H^2$.  For each integer $k\geq0$, our construction produces two 8-dimensional families of stable rank-2 co-Higgs bundles, one for $(c_1,c_2)=(0,k(k-1)H^2)$ and another for $(c_1,c_2)=(-H,k^2H^2)$.\\

\subsection{\fontfamily{cmss}\selectfont Facts about stability and bundles on \tb{CP}\textsuperscript2}

It will be useful to have at hand a couple of well-known facts about slope stability:\bitem\item For any stable $(V,\Phi)$, the subspace of $H^0(\End V)$ consisting of endomorphisms that commute with $\Phi$ is generated by $\mb 1_V$, and $(V,\Phi)$ is said to be \emph{simple}.  In particular, if $(V,\mb 0)$ is simple, then $V$ is said to be a \emph{simple vector bundle}.\item When $X$ is nonsingular and $V$ is a rank-2 bundle, we need only check \eqref{IneqStab} for sub-line bundles.
\eitem

\noin We will also use the following facts about bundles on $\CC\PP^2$: the only rank-2 locally-free sheaf with $(c_1,c_2)=(0,-H^2)$ is $\CO\plus\CO(-1)$; with $(c_1,c_2)=(0,0)$, only $\CO\plus\CO$; and with $(c_1,c_2)=(3H,3H^2)$, only $T$. 

\subsection{\fontfamily{cmss}\selectfont Some additional notation}

We will read $\End V\tensor T_X$ as $(\End V)\tensor T_X$.  We use $T$ without subscript to mean the tangent bundle of projective space.  If $V$ is a vector bundle on $\CC\PP^n$, then for economy we write $V(a)$ for $V\tensor\CO(a)$ and $\End V(a)$ for $(\End V)\tensor\CO(a)$.  We write $H^i(\CF)$ for the sheaf cohomology $H^i(X,\CF)$ whenever $X$ --- usually $\CC\PP^2$ or $\CC\PP^1\times\CC\PP^1$ --- is understood.  \\

\noin\tb{\fontfamily{cmss}\selectfont Acknowledgements.}  I thank Nigel Hitchin for introducing me to this topic and for his insights.  I acknowledge Nicolas Addington, Daniel Halpern-Leistner, Jonathan Fisher, Marco Gualtieri, Peter Gothen, Tam\'as Hausel, Lisa Jeffrey, Brent Pym, and Justin Sawon for useful discussions.  Thomas Peternell was indispensible in proving a vanishing theorem \citep[][Cor.9]{TP:11} required for a result in this paper.  Parts of this work were completed under funding from the Commonwealth Scholarship \& \hbox{Fellowship} Plan and the Natural Sciences \& \hbox{Engineering} Research Council of Canada.\\

\section{\fontfamily{cmss}\selectfont SIMPSON'S MODULI SPACES}

In \cite{TSCH:94I} and \cite{TSCH:94II}, Simpson constructs in two different ways a coarse moduli space of Higgs sheaves on a smooth projective variety of arbitrary dimension.  The first method is a direct GIT quotient giving the moduli space of coherent sheaves of $\Lambda$-modules on a projective variety $X$, where $\Lambda$ is a sheaf of $\CO_X$-algebras (possessing a filtration with certain properties). Taking $\Lambda=\tb{Sym}^\bullet(T_X)$ gives the moduli space of coherent Higgs sheaves on $X$.  In the second construction, Simpson passes to the \emph{spectral correspondence}: coherent Higgs sheaves on $X$ with fixed characteristic polynomial are identified with ordinary generically rank-1 coherent sheaves on a subvariety $S$ of a compactification of the cotangent bundle (supported away from the divisor at infinity).  Accordingly, the moduli space of Higgs sheaves on $X$ with fixed characteristic polynomial and the moduli space of sheaves of $\Lambda$-modules on $S$ for $\Lambda=\CO_S$ are isomorphic as varieties.  Allowing the characteristic polynomial to vary produces a fibration where the fibre is the moduli space of Higgs sheaves with a fixed characteristic polynomial and the base is the affine space of characteristic coefficients.  This is the \emph{Hitchin fibration}.  The total space of this fibration is isomorphic as a variety to the moduli space arising from Simpson's first construction.

For these moduli spaces, there is Gieseker's stability condition, which uses the Hilbert polynomial of a sheaf in lieu of the degree.  Although the examples we construct are slope stable, it is also well-known that slope stability implies Gieseker stability.  Therefore, if we replace $\Lambda=\tb{Sym}^\bullet(T_X)$ with $\Lambda=\tb{Sym}^\bullet(\Omega^1_X)$ in Simpson's first construction, or the cotangent bundle with the tangent bundle in the second construction, then the co-Higgs bundles constructed in this paper are points in one of Simpson's moduli spaces of $\Lambda$-modules.

Because of the great generality of Simpson's construction, it is difficult to extract concrete information about the global moduli space of co-Higgs bundles.  In particular, there is no obvious formula for the dimension of the space at a point in terms of the Chern classes of the underlying bundle.  Therefore, we focus on studying specific examples, and use deformation theory to see some of the local structure of the moduli space.

\begin{remark} The \emph{Hitchin map} on Simpson's moduli space --- the map sending a Higgs field to the coefficients of its characteristic polynomial --- is proper \citep[][Thm.6.11]{TSCH:94II}.  The families we construct will exhibit this property; however, it is important to note that because we are neither using Gieseker stability nor allowing sheaves that fail to be locally free, it cannot be expected \emph{a priori} for properness to be seen.\end{remark}

The subvariety $S$ of $T_X^\vee$ in Simpson's second construction is called a \emph{spectral cover}: it is a sheeted cover of $X$ whose sheets are eigenvalues of $\Phi$. In our case, the cover will be embedded in the total space of $T_X$.  We will use the following fact: if $\Phi\wedge\Phi=0$ and if the spectral cover belonging to $\Phi$ is smooth, then the dimension of $\ker[-,\Phi]$ is minimal at every point of $X$.  We will call such $\Phi$ \emph{regular}.  In particular, $\Phi_x$ may be nilpotent but $\Phi_x\neq0$ for any $x\in X$.  (See remarks in \cite{RD:95}, and in particular \citep[][Rmk.3.1]{KP:94}.)

While our focus is on Higgs fields taking values in the tangent bundle, it will be necessary in our arguments to consider Higgs fields taking values in a line bundle $L$.  If $\phi\in H^0(\Endz V\tensor L)$, then by $\phi^\vee$ we mean the dual element in the vector space $\Gamma(\Hom(V^\vee,V^\vee\tensor L))$.  Because we have $\Hom(V^\vee,V^\vee\tensor L)=\Hom(V,V\tensor L)$, we can regard $\phi$ and $\phi^\vee$ as elements of the same vector space, $H^0(\Endz V\tensor L)$.

\begin{lemma}\label{LemmRegular} Let $(V,\phi)$ be a regular Higgs bundle on a smooth complex manifold with $\phi\in H^0(\emph{\End}_0V\tensor L)$. Then, we have a short exact sequence\beqn0\longrightarrow L^\vee\stackrel{\phi}{\longrightarrow}\emph{\End}_0V\stackrel{[-,\phi]}{\longrightarrow}Q\longrightarrow 0,\eeqn\noin where $Q$ is the sheaf-theoretic image of $[-,\phi]$ in $\emph{\End}_0V\tensor L$.  Now regard $\phi$ as an element of $\Gamma(\emph{\mbox{Hom}}(\emph{\End}_0V\tensor L,L^{\tensor2}))$; if for some $c\in\CC$ we have $\phi^\vee=c\phi$, then $\emph{\ker}\phi\cong Q$.\end{lemma}

\begin{proof} Begin by considering $\phi$ as an element of $\Gamma(\mbox{Hom}(L^\vee,\Endz V))$.  The image $\phi(L^\vee)$ is in the kernel of $[-,\phi]:\Endz V\rightarrow\Endz V\tensor L$.  The regularity of $\phi$ means that $\dim[-,\phi]$ is minimal, and so $\img\phi=\ker[-,\phi]$, producing for us the short exact sequence.  Now assume that $\phi^\vee=c\phi$. Take the maps $L^\vee\stackrel{\phi}{\longrightarrow}{\Endz}V\stackrel{[-,\phi]}{\longrightarrow}\Endz V\tensor L$ and complete them to a four-term exact sequence\beqn0\longrightarrow L^\vee\stackrel{\phi}{\longrightarrow}{\Endz}V\stackrel{[-,\phi]}{\longrightarrow}\Endz V\tensor L\longrightarrow M\longrightarrow0,\label{PreTwist}\eeqn\noin wherein $M$ is a line bundle.  The dual of this sequence is\beqn0\longrightarrow M^\vee\stackrel{\phi^\vee}{\longrightarrow}{\Endz}V\tensor L^\vee\stackrel{[\phi^\vee,-]}{\longrightarrow}\Endz V\longrightarrow L\longrightarrow0,\nonumber\eeqn\noin which is equivalently\beqn0\longrightarrow M^\vee\stackrel{c\phi}{\longrightarrow}{\Endz}V\tensor L^\vee\stackrel{-c[-,\phi]}{\longrightarrow}\Endz V\longrightarrow L\longrightarrow0.\label{PostTwist}\eeqn\noin  The sequences \eqref{PreTwist} and \eqref{PostTwist} differ only by a twist by $L^\vee$, meaning that $M=L^{\tensor 2}$.  Then \eqref{PreTwist} becomes \beqn0\longrightarrow L^\vee\stackrel{\phi}{\longrightarrow}{\Endz}V\stackrel{[-,\phi]}{\longrightarrow}\Endz V\tensor L\stackrel{\phi}{\longrightarrow}L^2\longrightarrow0\nonumber\eeqn\noin from which the second claim in the statement of the lemma must follow.\end{proof} 

\section{\fontfamily{cmss}\selectfont DEFORMATIONS OF CO-HIGGS BUNDLES ON COMPLEX\\MANIFOLDS}\label{Defs}

Let $X$ be a complex manifold of any dimension; $(V,\Phi)$, a co-Higgs bundle over $X$.  The condition $\Phi\wedge\Phi=0$ makes $\wedge\Phi$ into a differential on \v Cech cochains for the bundle $\Endz  V\tensor\wedge^{\bullet}T_X$.  The operation $\wedge\Phi$ commutes with the \v Cech coboundary $\delta$, making the total module\beqn\ds(\tb{C}^{\bullet}(\Endz  V\tensor\wedge^{\bullet}T_X);\;D=\delta+\wedge\Phi)\nonumber\eeqn\noin into a first-quadrant double complex.

A spectral sequence is defined by choosing the 0-th page to be the module\beqn\ds(\CE_0^{p,q}=\tb{C}^{\bullet}(\Endz  V\tensor\wedge^{\bullet}T_X);\;d_0=\delta),\nonumber\eeqn\noin noting that $d_0:\CE_0^{p,q}\rightarrow\CE_0^{p,q+1}$.  Then proceed by setting:

\bitem
\item[\tb{(i)}] $\CE_1^{p,q}=H^q_{d_0}(\CE_0^{p,\bullet})=H^q(\Endz  V\tensor\wedge^pT_X)$
\item[\tb{(ii)}] $d_1=\wedge\Phi:\CE_1^{p,q}\rightarrow\CE_1^{p+1,q}$ (acting on \v Cech $q$-cochains)
\item[\tb{(iii)}] $\ds\CE_2^{p,q}=H^p_{d_1}(\CE_1^{\bullet,q})=\frac{\ker H^q(\Endz  V\tensor\wedge^pT_X)\stackrel{\wedge\Phi}{\longrightarrow}H^q(\Endz  V\tensor\wedge^{p+1}T_X)}{\img H^q(\Endz  V\tensor\wedge^{p-1}T_X)\stackrel{\wedge\Phi}{\longrightarrow}H^q(\Endz  V\tensor\wedge^pT_X)}$
\eitem

In addition to $p,q\geq0$, the sequence enjoys these finiteness properties: $H^q(\Endz  V\tensor\wedge^pT_X)=0$ when either $p,q>\dim(X)$.  Note that $d_2:\CE_2^{p,q}\rightarrow\CE_2^{p+2,q-1}$ is given by\beqn d_2(\psi) & = & \theta\wedge\Phi\label{d2},\eeqn\noin where $\theta\in\tb{C}^{q-1}(\Endz  E\tensor\wedge^{p+1}T)$ is the solution of the equation $\psi\wedge\Phi-\delta\theta=0\in\tb{C}^q(\Endz  E\tensor\wedge^{p+1}T_X)$.  We are interested in $d_2$ because of the following fact from homological algebra:

\begin{proposition}\label{PropFiveTerm}  If $\HH^\bullet$ is the hypercohomology of the double complex, then there is an exact sequence\beqn0\longrightarrow\CE_2^{1,0}\longrightarrow\HH^1\longrightarrow\CE_2^{0,1}\stackrel{d_2}{\longrightarrow}\CE_2^{2,0}\longrightarrow\HH^2\label{ESFiveTerm}.\eeqn\end{proposition}

For the applications we have in mind, we will always have $d_2|_{\CE_2^{0,1}}=0$.  Note that, in this case, a first-order deformation of $(V,\Phi)$ has two components: a deformation in $\CE_2^{1,0}$ and a deformation in $\CE_2^{0,1}$.  First-order deformations of the Higgs field that are holomorphic with respect to the given complex structure on $V$ are given by elements of\beqn\CE_2^{1,0} & = & \frac{\ker H^0(\Endz  V\tensor T_X)\stackrel{\wedge\Phi}{\longrightarrow}H^0(\Endz  V\tensor\wedge^2T_X)}{\img H^0(\Endz  V)\stackrel{\wedge\Phi}{\longrightarrow}H^0(\Endz  V\tensor T_X)}.\label{DefnE10}\eeqn\noin  On the other hand,\beqn\CE_2^{0,1} & = & \ker H^1(\Endz  V)\stackrel{\wedge\Phi}{\longrightarrow}H^1(\Endz  V\tensor T_X)\label{DefnE01}\eeqn\noin is a space of Kodaira-Spencer classes for $V$, but only those corresponding to first-order deformations of the bundle $V$ along which the given $\Phi$ remains holomorphic.

\begin{remark}The deformation theory for Higgs bundles on nonsingular curves appears in the works of Nitsure, Biswas and Ramanan, and Bottacin; respectively \cite{NN:91}, \cite{BR:94}, and \cite{FB:95}.  The theory for Higgs bundles on curves can be recovered from the sequences above by replacing $T_X$ with $T_X^\vee=\w_X$ throughout, and noting that nonzero terms on the $\CE_2^{\bullet,\bullet}$ page will be concentrated in the band $0\leq p,q\leq 1$.\end{remark}

\section{\fontfamily{cmss}\selectfont VANISHING THEOREMS}

On surfaces, the presence of stable co-Higgs bundles is skewed to the nonpositive end of the Kodaira spectrum, at least for rank 2.  Theorem \ref{ThmVanishing} below supports this.

\begin{lemma}\label{LemmDet0} Let $X$ be a nonsingular complex projective surface with a rank-2 co-Higgs bundle $(V,\Phi)\rightarrow X$ for which $\det(\Phi)=0$.  If $\Phi$ is not identically zero, then there exist line bundles $L$ and $M$ on $X$ with the following properties: $V$ is an extension\beqn0\rightarrow L\rightarrow V\rightarrow M\tensor\MF I_Z\rightarrow0\nonumber\eeqn\noin in which $\MF I_Z$ is an ideal sheaf of points $Z\subset X$; $L=\emph{\ker}\Phi$; and $\Phi$ is a global holomorphic section of $M^\vee\tensor L\tensor T_X$.\end{lemma}

\begin{proof}  Since $V$ has rank 2, $\mbox{tr}(\Phi)=0$ and $\det(\Phi)=0$ together imply that $\Phi$ is nilpotent. Since $\Phi$ itself is not identically zero, there must exist a line bundle $L=\ker\Phi$ included as a sheaf in $V$, and therefore a short exact sequence of sheaves\beqn0\rightarrow L\rightarrow V\rightarrow M\tensor\MF I_Z\rightarrow0\nonumber\eeqn\noin in which $M$ is a line bundle on $X$ and $Z\subset X$ is a set of points.  It follows that\beqn\Phi\in H^0(X,(M\tensor\MF I_Z)^\vee\tensor L\tensor T_X).\nonumber\eeqn\noin  We may extend $\Phi$ uniquely over $Z$ by the theorem of Hartogs, and so we have \hbox{$\Phi\in H^0(X,M^\vee\tensor L\tensor T_X)$,} once we agree to reuse $\Phi$ for the extension.\end{proof}

\begin{theorem}\label{ThmVanishing} Let $\iota:X\inject\CC\PP^N$ be a nonsingular, connected surface of general type, and let $(V,\Phi)\rightarrow X$ be a semistable rank-2 co-Higgs bundle with $c_1(V)=0$ or $c_1(V)=-H$, where $H=c_1(\iota^*\CO(1))$.  Then, $\Phi=0$.\end{theorem}

\begin{proof} We do $c_1(V)=0$ first.  With $X$ as in the statement, we must have $H^0(X,\mb{S}^2T_X)=0$.  This follows from a more general vanishing result of Peternell \citep[][Cor.9]{TP:11}, saying that $H^0(X,T^{\tensor m})=0$ for all $m\geq1$ when $X$ is of general type.  This means that $\det\Phi=0$.  By Lemma \ref{LemmDet0}, $V$ has a sub-line bundle $L$ and $\Phi$ is in $H^0(X,L^{2}\tensor T_X)$.  Since $(V,\Phi)$ is semistable, it must follow that $\deg L\leq0$, which in turn means $\deg L^{-2}\geq0$.  By definition, this means that $L^{-2}.C\geq0$ for any curve $C$ in the linear system $\abs{\iota^*\CO(1)}$, and so it follows that $L^{-2}$ is pseudo-effective (see Theorem 0.2 and Corollary 0.3 of \cite{SB:13}).   Peternell shows that $H^0(X,T_X\tensor D^\vee)=0$ for any pseudo-effective line bundle $D$ on a projective manifold of general type \citep[][Cor.9]{TP:11}.  If we take $D=L^{-2}$, then $H^0(L^2\tensor T_X)=0$, and $\Phi$ must vanish identically.

In the $c_1(V)=-H$ case, we have $\Phi\in H^0(X,M^\vee\tensor L\tensor T_X)$ where $c_1(M)=-(1+k)H$ and $c_1(L)=kH$ for some $k$, and so $\deg(M^\vee L)=1+2k$.  For semistability, we need $\deg L=k\leq-1/2$, i.e. $k\leq-1$.  This means that $\deg(M^\vee L)^\vee\geq1$, and the remainder of the argument proceeds as in the even case.\end{proof}

A similar theorem holds for K3 surfaces, but the proof has a different flavour.  We will give an outline, in the case of $c_1(V)=0$.  When $X$ is K3, we have $H^0(\tb{S}^2T_X)=0$. From Lemma \ref{LemmDet0}, if $X$ admits a stable rank-2 co-Higgs bundle $(V,\Phi)$, then there exists a line bundle $L$ included as a subsheaf in $V$, and $\Phi$ is an element of $H^0(X,L^2\tensor T_X)$.  We can use the vanishing theorem of Kobayashi and Wu \citep[][p.1]{KW:70} to show $H^0(L^2\tensor T_X)=0$, by proving that there exists a certain curvature $(1,1)$-form $F$ on $L$ that is negative definite after being contracted with the metric coming from the K\"ahler form $\w$ on $X$.  (The vanishing of the Ricci tensor means that the curvature on $L^{2}\tensor T_X$ comes from $L^{2}$.)  We can construct $F$ out of any curvature form $F_0$ on $L$, using the fact that the stability condition $\deg(L)<0$ means that every $F_0$ must satisfy $\int_XF_0\wedge\w=\int_X[c_1(L)].[\w]=c<0$.  We use Hodge theory to produce a function $h$ such that $F=F_0+\di\bar\di h$ is a tensor with the desired properties.

A reason to posit that $\CC\PP^2$ might be a generous source of co-Higgs bundles is that the vanishing theorem ties the existence of stable rank-2 examples to the availability of holomorphic sections of $\tb{S}^2T_X$. The projective plane has many.

\section{\fontfamily{cmss}\selectfont THREE EXAMPLES}

From now on, $\mb{CHB}(k_1,k_2)$ stands for the moduli space of stable rank-2 co-Higgs bundles on $\CC\PP^2$ with Chern classes $(c_1,c_2)=(k_1H,k_2H^2)$.

\subsection{\fontfamily{cmss}\selectfont Examples 1 and 2: decomposable cases}

It is natural to start with extensions of one line bundle by another, say, $L_1$ by $L_2$.  On $\CC\PP^2$, the only such extensions are the trivial ones. Not every direct sum, however, admits a stable $\Phi$.

\begin{proposition}\label{PropDirect}  Suppose that there exists a stable $\Phi\in H^0(\CC\PP^2,\emph{\mbox{End}}_0V\tensor T)$ for $V=\CO(m_1)\plus\CO(m_2)$.  Then $\abs{m_1-m_2}\leq1$.\end{proposition}

\begin{proof}Consider the Euler sequence on $\CC\PP^2$:\beqn0\longrightarrow\CO\longrightarrow\Plus_{i=1}^3\CO(1)\longrightarrow T\longrightarrow0\label{SESEuler}.\eeqn\noin  If we twist the terms of the sequence by $\CO(-d)$ for any $d>1$, then the free terms become $\CO(-d)$ and $\Plus_{i=1}^3\CO(1-d)$, respectively, which are sums of negative-degree line bundles only.  Therefore, $T(-d)$ has no global sections for $d>1$.

Assume without loss of generality that $m_1\geq m_2$.  The Higgs field $\Phi$ has a component $\psi:\CO(m_1)\rightarrow T(m_2)\in H^0(T(m_2-m_1))$.  If $m_1-m_2>1$, then $\psi=0$ and $\CO(m_1)$ is invariant and destabilizing, contradicting the stability of $(V,\Phi)$.\end{proof}

After we impose $c_1(V)=-H$ or $c_1(V)=0$, stability permits only $V=\CO\plus\CO(-1)$ or $V=\CO\plus\CO$, respectively.

We begin with $V=\CO\plus\CO(-1)$. Each $\Phi\in H^0(\Endz V\tensor T)$ takes the form\begin{center}\bmath\Phi\;\;=\;\;\left(\begin{array}{ccc}A & & \;\;B\\C & & -A\end{array}\right)\emath\end{center}\noin for some $A\in H^0(T)$, $B\in H^0(T(1))$, and $C\in H^0(T(-1))$.  This is a stable Higgs field for $V$ if and only if $C$ is not identically zero, so that the trivial sub-line bundle in $V$ is not preserved.  The pair $(V,\Phi)$ is a stable co-Higgs bundle if and only if $C\neq0$ and the form\begin{center}\bmath\Phi\wedge\Phi\;\;=\;\;\left(\begin{array}{cc}B\wedge C & \;\;2A\wedge B\\2C\wedge A & \;\;C\wedge B\end{array}\right)\emath\end{center}\noin vanishes identically.  This vanishing is equivalent to $A$, $B$, and $C$ satisfying the simultaneous system\beqn A\wedge B\;\;=\;\;0,\;\;\;\;A\wedge C\;\;=\;\;0,\;\;\;\;B\wedge C\;\;=\;\;0\,.\nonumber\eeqn\noin  Since $C$ is not identically zero, $C$ vanishes on a single point $p\in\CC\PP^2$. Away from $p$, the simultaneous conditions imply that $A=\lambda C$ and $B=\mu C$, where $\lambda$ is a section of $\CO(1)$ and $\mu$ is a section of $\CO(2)$ over $\CC\PP^2\bsh\set{p}$.  Hartogs' theorem allows us to extend each of $\lambda$ and $\mu$ uniquely to sections over the whole of $\CC\PP^2$.  Thus, every stable $\Phi$ satisfying $\Phi\wedge\Phi=0$ can be written\beqn\Phi\;\;=\;\;\phi\tensor C\;\;=\;\;\left(\begin{array}{cc}\lambda & \;\;\mu\\1 & -\lambda\end{array}\right)\tensor C,\nonumber\eeqn\noin where $C\in H^0(T(-1))\bsh\set{0}$ and the matrix part is a section $\phi\in H^0(\Endz V(1))$.  Using the automorphism\beqn\Psi & = & \left(\begin{array}{cc} 1 & \;\;\;\;\lambda\\0 & \;\;\;\;1\end{array}\right)\nonumber\eeqn\noin of $V=\CO\plus\CO(-1)$, we can transform $\Phi$ within its equivalence class to\beqn\Psi^{-1}\Phi\Psi & = & \left(\begin{array}{cc}0 & \;\;\;\;q\\1 & \;\;\;\;0\end{array}\right)\tensor C\;\;=\;\;\left(\begin{array}{cc}0 & \;\;\;\;qC\\C & \;\;\;\;0\end{array}\right),\nonumber\eeqn\noin where $q=-\det\phi=\lambda^2+\mu\in H^0(\CO(2))\cong\CC^6$.  It is clear that the data $(q,C)$ determines a unique Higgs field $\Phi$, but not vice-versa.  If we scale $C$ by any $t\in\CC^*$ and $q$ by $t^{-2}$, then we obtain the same $\Phi$.  (Equivalently, the family of automorphisms $\Psi_t=\mbox{diag}(t,t^{-1})$ fixes $\Phi$.)   In other words, the moduli space of stable, integrable Higgs fields for $V=\CO\plus\CO(-1)$ is a quotient of $\CC^6\times\CC^3$ by $\CC^*$ acting with weight $-2$ on $\CC^6$ and with weight $1$ on $\CC^3$.  According to the stability condition $C\neq0$, this quotient is $(\CC^6\times\CC^3\bsh\set{0})//\CC^*$, which is isomorphic to the total space of the rank-6 vector bundle $\CO_{\CC\PP^2}(-2)^{\plus6}$.  (The opposite linearization, $\CC^3\bsh\set{0}\times\CC^3$, would have resulted in an orbi-line bundle.  This quotient problem is discussed in $\S$2.4 of \cite{RT:11}.)  We package this discussion as

\begin{theorem}\label{Thmk=0}  The moduli space $\emph{\tb{CHB}}(-1,0)$ is isomorphic to the total space of $\CO_{\CC\PP^2}(-2)^{\plus6}$.\end{theorem}

The $V=\CO\plus\CO$ case has a complication.  Now, a Higgs field $\Phi$ is a matrix\begin{center}\bmath\Phi\;\;=\;\;\left(\begin{array}{cc}A & \;\;\;\;B\\C & \;\;-A\end{array}\right)\emath\end{center}\noin whose entries $A$, $B$, $C$ are holomorphic vector fields.  There are no unstable Higgs fields whatsoever, but there are semistable ones that are not stable --- in particular, $\Phi=0$.  More generally, when one of $B$ or $C$ is identically zero, a degree 0 sub-line bundle will be preserved.  Recall that slope-semistable objects are subject to \emph{S-equivalence}, first introduced for vector bundles in \cite{CSS:67}, identifying those points whose associated graded objects are identical. (The underlying bundle $V=\CO\plus\CO$ is fixed, so we need only concern ourselves with identifying associated graded Higgs fields.)

Let\beqn\widetilde\CS & = & \left.\set{\left.\Phi_{q,C}=\left(\begin{array}{ccc} 0 & ~~ & q\\ 1 & ~~ & 0\end{array}\right)C\;\right|\;q\in\CC\mbox{ and }C\in H^0(T)\bsh\set{0}}\,\right/\,\sim,\nonumber\eeqn\noin where $\sim$ is the $\CC^*$ action defined by letting matrices with the form $\psi_t=\left(\begin{array}{ccc} t & ~~ & 0 \\ 0 & ~~ & t^{-1}\end{array}\right)$ act by conjugation on the $\Phi_{q,C}$.  The action of $\Psi_t=\left(\begin{array}{ccc} t & ~~ & 0 \\ 0 & ~~ & t^{-1}\end{array}\right)$ identifies $\Phi_{q,C}$ with $\Phi_{t^{-4}q,t^2C}$.  In other words, $\CC^*$ acts with weight $-2$ on $q$ and $+1$ on $C$, and we have $\widetilde\CS\cong\mbox{Tot}(\CO_{\CC\PP^2}(-2)^{\plus6})$, cf. \citep[][$\S$2.4]{RT:11}.  Let $\widetilde\CS_0$ stand for the contraction of the zero section to a point, which we represent by $\Phi=0$.

\begin{theorem}\label{Thmk=1}   $\widetilde\CS$ is an open dense subset of $\mb{CHB}(0,0)$.\end{theorem}

\begin{proof} Let $\set{s=0}$ stand for the zero section of $\widetilde\CS\cong\mbox{Tot}(\CO_{\CC\PP^2}(-2)^{\plus6})$. It is clear that $\widetilde\CS\bsh\set{s=0}$ is a subset of $\mb{CHB}(0,0)$, for after quotienting $\widetilde S$ by $\sim$, the assignment of an element $[\Phi_{q,C}]\in\widetilde S_0$ to its determinant $q\tensor C\tensor C\in H^0(\mb S^2T)$ is injective.  Along the zero section, the determinant is constant (and equal to $0$), but the corresponding Higgs fields $\Phi_{0,C}$ are non-isomorphic for different $C\in\CC\PP^2$.  However, S-equivalence replaces all of the $\Phi_{0,C}$ with the zero Higgs field. Therefore, while $\widetilde\CS$ is not a subset of $\mb{CHB}(0,0)$, the contraction $\widetilde\CS_0$ is.

  Now, consider those $\Phi=\left(\begin{array}{ccr} A & ~~ & B\\ C & ~~ & -A\end{array}\right)$ for which $C$ is not identically zero and vanishes at a single point in $\CC\PP^2$.  We use $\CS$ to denote the set of such $\Phi$ in $H^0(\Endz V\tensor T)$.
By the same argument as for for $\CO\plus\CO(-1)$, the solutions of $\Phi\wedge\Phi=0$ in $\CS$ are those $A,B,C$ for which $A=aC$ and $B=bC$, where $a,b\in\CC$.
A gauge transformation $\Psi=\left(\begin{array}{ccr} 1 & ~~ & a\\ 0 & ~~ & 1\end{array}\right)$ takes a solution $\Phi$ to $\Phi_{q,C}=\left(\begin{array}{ccc}0 & ~~ & q\\ 1 & ~~ & 0\end{array}\right)C$, where $q=a^2+b\in\CC$.  Again, the isomorphism class of $\Phi$ does not determine $(q,C)$ uniquely.  To remedy this, we take the quotient $\CS/\sim$, and then identify those points of the form $\Phi_{0,C}$, as per S-equivalence.  The final quotient is a proper subset of $\widetilde\CS_0$.  
That $\widetilde\CS_0$ is open dense in $\mb{CHB}(0,0)$ comes from the fact that $C\neq0$ is generic.   \end{proof} 

\begin{corollary}\label{Cork=1} The moduli space $\mb{CHB}(0,0)$ is $8$-dimensional.\end{corollary}

\begin{corollary} The moduli spaces $\mb{CHB}(-1,0)$ and $\mb{CHB}(0,0)$ are not isomorphic as varieties.\end{corollary}

\begin{remark}  Because of the zero Higgs field, $\widetilde\CS_0$ is not a subvariety of a Simpson moduli space.  (Slope stable implies Gieseker stable, but slope semistable does not imply Gieseker semistable in general.)  The set $\widetilde\CS$ with the zero section excised completely is, on the other hand, contained in a Simpson moduli space.\end{remark}

\subsection{\fontfamily{cmss}\selectfont Example 3: the tangent bundle}

After direct sums of two line bundles, the natural rank-2 vector bundle to consider is the tangent bundle itself, which for $\CC\PP^2$ is indecomposable.  Unlike the direct sums, there is no stability condition to solve: $T$ is stable as a vector bundle, and therefore $(T,\Phi)$ is stable for any $\Phi$.  Note that $H^0(\CC\PP^2,\Endz T\tensor T^*)=0$: the tangent bundle of $\CC\PP^2$ fails to admit any nonzero Higgs fields in the conventional sense.  On the other hand, the vector space $H^0(\Endz T\tensor T)$ is 18-dimensional.

The space $H^0(\Endz T(1))$ is 6-dimensional, and we have a canonical isomorphism\beqn\CC^3\tensor H^0(\Endz T(1)) & \cong & H^0(\Endz T\tensor T).\label{IsoO(1)ToTVer2}\eeqn\noin This isomorphism comes to us by way of the Euler sequence \eqref{SESEuler}.  Applying $\Endz T\tensor$ to \eqref{SESEuler} produces another short exact sequence,\beqn0\longrightarrow\Endz T\longrightarrow(\Endz T(1))^{\plus 3}\longrightarrow\Endz T\tensor T\longrightarrow0.\nonumber\eeqn\noin The first four terms in cohomology are\beqn 0\rightarrow H^0(\Endz T)\rightarrow H^0(\Endz T(1))^{\plus 3}\rightarrow H^0(\Endz T\tensor T)\rightarrow H^1(\Endz T).\nonumber\eeqn\noin The leftmost $H^0$ is $\set{0}$ because $T$ is stable.  The space $H^1(\Endz T)$ is $\set{0}$, as $T$ is rigid.  What remains is the isomorphism \eqref{IsoO(1)ToTVer2}. It also follows from the Euler sequence that the $\CC^3$ in \eqref{IsoO(1)ToTVer2} is identified with $H^0(T(-1))$.

If $\phi$ is any element of $H^0(\Endz T(1))$, then $(T,\phi)$ is a stable $O(1)$-valued Higgs bundle.  If $\det\phi\in H^0(\CC\PP^2,\CO(2))$ is irreducible as a polynomial, then the characteristic equation of $\phi$ determines a nonsingular spectral cover of $\CC\PP^2$ embedded in the total space of $\CO(1)$.  Irreducibility of $\det\phi$ is an open condition, and so the generic $\phi$ is regular.

According to \eqref{IsoO(1)ToTVer2}, if $\set{\phi_1,\dots,\phi_6}$, $\set{C_1,C_2,C_3}$ are bases for $H^0(\Endz T(1))$ and $H^0(T(-1))$, respectively, then any $\Phi\in H^0(\Endz T\tensor T)$ can be written as a tensor $\Phi=\sum_{i=1}^6\sum_{j=1}^3a_{ij}\phi_i\tensor C_j$ for some $a_{ij}\in\CC$. In particular, we can take a basis of $H^0(\Endz T(1))$ consisting of regular elements $\phi_1,\dots,\phi_6$.  Consider Lemma \ref{LemmRegular} applied to $\phi_i$: there is a short exact sequence\beqn0\longrightarrow\CO(-1)\stackrel{\phi_i}{\longrightarrow}\Endz T\stackrel{[-,\phi_i]}{\longrightarrow}Q\rightarrow 0,\nonumber\eeqn\noin where $Q=\img[-,\phi_i]\subset\Endz T(1)$.  Twisting by $\CO(1)$ gives an equivalent sequence\beqn0\longrightarrow\CO\stackrel{\phi_i}{\longrightarrow}\Endz T(1)\stackrel{[-,\phi_i]}{\longrightarrow}Q(1)\rightarrow 0,\nonumber\eeqn\noin where $Q(1)$ is now the image of $[-,\phi_i]$ when it is regarded as a map of sheaves from $\Endz T(1)$ to $\Endz T(2)$.  In cohomology, we have\beqn0\longrightarrow\CC\stackrel{\phi_i}{\longrightarrow}H^0(\Endz T(1))\stackrel{[-,\phi_i]}{\longrightarrow}H^0(Q)\longrightarrow0.\nonumber\eeqn\noin  The kernel in the sequence consists of the scalar multiples of $\phi_i$.  Consequently, $H^0(Q)=\img[-,\phi_i]\subset H^0(\Endz T(2))$ is 5-dimensional and spanned by the $[\phi_j,\phi_i]$, $j\neq i$. 

Therefore, $\Phi\wedge\Phi$ is a linear combination of terms $[\phi_i,\phi_j]C_k\wedge C_l$, with $i\neq j$ and $k\neq l$.  In order to have $\Phi\wedge\Phi=0$, we must have either $\Phi=\sum a_j\phi_i C_j$ for a fixed $i$ or $\Phi=\sum a_i\phi_i C_j$ for a fixed $j$.  In either case, $\Phi=\phi\tensor C$ for some $\phi\in H^0(\Endz T(1))\cong\CC^6$ and for some $C\in H^0(T(-1))\cong\CC^3$.  

As with the $\CO\plus\CO(-1)$ case, the moduli space is a quotient of $\CC^6\times\CC^3$, but now the $\CC^6$ cannot be identified with the space of determinants for $\phi$.  In the case of $V=T$, we have that $\det:H^0(\Endz T(1))\rightarrow H^0(\CO(2))$ is a double cover of $\CC^6$ by itself.  Of course, this is true in the case of $V=\CO\plus\CO(-1)$, but for that bundle we can use the automorphism $\psi=\mbox{diag}(1,-1)$ to identify $\phi$ and $-\phi$ in the moduli space.  Consequently, conjugacy classes of $\CO(1)$-valued Higgs fields for $\CO\plus\CO(-1)$ are in 1:1 correspondence with their determinants.  In contrast, there are no automorphisms of $T$ other than multiples of $\mb 1_T$ (by the \emph{stable implies simple} property), and so $\det\phi$ does not determine the isomorphism class of $\phi$, unless $\det\phi=0$.

The result is that we obtain the moduli space by quotienting \beqn H^0(\Endz T(1))\times H^0(T(-1)) & = & \CC^6\times\CC^3\nonumber\eeqn\noin by a $\CC^*$ action with weights $\pm1$, which accounts for the fact that $(\lambda^{-1}\phi)\tensor(\lambda C)$ gives rise to the same $\Phi$ as $\phi\tensor C$. As before, we need a stability condition for the quotient, but this time it does not descend automatically from the slope stability condition.  Both $C\neq0$ and $\det\phi\neq0$ are acceptable conditions, and neither contradicts slope stability.  Staying consistent with the previous examples, we take $C\neq0$, resulting in the bundle $\CO_{\CC\PP^2}(-1)^{\plus6}$.  As in the $V=\CO\plus\CO$ case, we need to contract the zero section to a point --- however, the reason is different.  Previously, there existed points that were semistable but not stable.  There are no such points for $V=T$.  The problem is that the points along the zero section are of the form $\phi\tensor C$ with $\det\phi=0$.  Since $\det\phi=0$ implies $\phi=0$, every point on the zero section must be the zero Higgs field.

\begin{theorem}\label{Thmk=2} The moduli space $\mb{CHB}(-1,1)$ is the total space of $\CO_{\CC\PP^2}(-1)^{\plus6}$ with points along the zero section identified with $\Phi=0$.\end{theorem}

\begin{remark} The moduli space we have just described is $\mb{CHB}(3,3)$: we found all of the integrable Higgs fields for $T$, and $T$ is the only rank-2 bundle with Chern classes $c_1=3H$, $c_2=3H^2$.  The normalizing isomorphism $\mb{CHB}(3,3)\cong\mb{CHB}(-1,1)$ comes from tensoring $T$ by $\CO(-2)$.  Tensoring by $\CO(-3)$ instead, we get the moduli space $\mb{CHB}(-3,3)$ of co-Higgs bundles with underlying bundle isomorphic to the cotangent bundle.\end{remark}

\section{\fontfamily{cmss}\selectfont SCHWARZENBERGER BUNDLES}

There is a common outcome in the examples so far: each of $\CO\plus\CO(-1)$, $\CO\plus\CO$, and $T$ underlies an $8$-dimensional family of co-Higgs bundles.  For two of these examples, stable integrable Higgs fields always decompose as $\Phi=\phi\tensor C$, with $\phi$ an $\CO(1)$-valued Higgs field and $C$ a section of $T(-1)$.  It turns out that there is a framework into which these examples can be placed, one that provides (a) a rationale for the decomposition of $\Phi$, and (b) many more examples.

We recall basic facts surrounding Schwarzenberger's construction of rank-2 holomorphic vector bundles on $\CC\PP^2$  \cite{RLES:61I}.  For each nonzero irreducible polynomial $\p\in H^0(\CC\PP^2,\CO(2))$, we can find a holomorphic cover\beqn f^\p\,:\,\CC\PP^1\times\CC\PP^1\stackrel{2:1}{\longrightarrow}\CC\PP^2,\nonumber\eeqn\noin branched over the nonsingular conic determined by $\p$.  (Two $\p$ that differ only in scale determine the same cover and the same branch conic.)  Each $\p$ also comes with a sequence of sheaves $\set{V^\p_k}_{k\ge0}$, defined by\beqn V^\p_k:=f^\p_*\CO(0,k)\nonumber\eeqn\noin for $k\geq0$, where $\CO(0,k):=\CO_1\tensor\CO_2(k)$ has $\CO_1$ is pulled back from one ruling and $\CO_2$ from the other.  The sheaves $V^\p_k$ are now called \emph{Schwarzenberger bundles}.  To emphasize the choice of $k$, we may refer to $V_k^\p$ as a ``$k$-Schwarzenberger'' bundle.  For convenience, we will use the symbol $[\p]$ to refer to both the projective class of a section $\p\in H^0(\CC\PP^2,\CO(2))$ and to the geometric conic $\set{x\in\CC\PP^2\;:\;\p(x)=0}$.

The following properties of $V^\p_k$ can be found in \cite{RLES:61I}, \citep[][pp.46--51]{RF:98}, and \citep[][\S2]{BS:92}:

\bitem
\item[\tb{(i)}] The naming is sound: $V^\p_k$ is locally-free of rank 2 for all $[\p]\in\CC\PP^5=\PP\left(H^0(\CO(2))\bsh\set{0}\right)$ and for all $k\geq0$ \citep[][Thm.2]{RLES:61I}.
\item[\tb{(ii)}] For $k=0,1,2$, the bundle $V^\p_k$ is rigid.  For $k\geq$, the space of first-order deformations of $V^\p_k$ is\beqn H^1(\CC\PP^2,\Endz V^\p_k) & = & \CC^{k^2-4}.\nonumber\eeqn\noin (These results will be recovered in our calculations below, in $\S$6.1.)
\item[\tb{(iii)}] $V^\p_k$ is indecomposable and slope stable for $k\geq2$ \citep[][Thm.2.7]{BS:92}.
\item[\tb{(iv)}] $V^\p_0\cong\CO\plus\CO(-1)$, $V^\p_1\cong\CO\plus\CO$, and $V^\p_2\cong T(-1)$ are independent of $[\p]$, even when $[\p]$ is reducible.  (This is true because $\CO\plus\CO(-1)$, $\CO\plus\CO$, and $T$ are the only locally-free sheaves on $\CC\PP^2$ with their Chern classes.  In the case where $[\p]$ is reducible, the quadric double cover is singular and the line bundle $\CO(0,k)$ is replaced by a reflexive sheaf $I$ with $c_1(I)=c_1(\CO(0,k))$.  The direct image of $I$ is reflexive, and hence locally free on $\CC\PP^2$.)
\item[\tb{(v)}] For $k\geq3$, $V^\p_k\cong V^{\p'}_k$ if and only if $V^\p_k$ and $V^{\p'}_k$ come from the same branch conic, that is, if and only if $\p\equiv\p'$ in $\CC\PP^5$.  (The line bundle $\CO(0,k)$ is rigid, and so the only data that goes into constructing $V^\p_k$ is $[\p]\in\CC\PP^5$.  On the other hand, for $k\geq3$, $h^1(\Endz V^\p_k)=k^2-4\geq5$.)
\item[\tb{(vi)}] After a Grothendieck-Riemann-Roch calculation \citep[][Thm.5]{RLES:61I}, one obtains\beqn c_1(V^\p_k) & = & (k-1)H\nonumber\\c_2(V^\p_k) & = & \frac{k(k-1)}{2}H^2\nonumber\eeqn\noin and\beqn c_1((V^\p_k)^\vee) & = & (1-k)H\nonumber\\c_2((V^\p_k)^\vee) & = & \frac{k(k-1)}{2}H^2.\nonumber\eeqn\noin
\item[\tb{(vii)}] $V^\p_k\cong V^\p_{k'}$ if and only if $k=k'$.  This follows from the Chern data: even after normalizing $c_1$ to one of $0$ or $-H$, $c_2$ remains a strictly monotone function of $k$.
\eitem

Note that for $k>3$, $h^1(\Endz V^\p_k)>5$, and so there are deformations of $V^\p_k$ that are not obtained from the Schwarzenberger construction; in other words, there are deformations $V$ of $V^\p_k$ for which $V\ncong V^{\p'}_k$ for any $[\p']\in\CC\PP^5$.

Having studied $k=0,1,2$ already, we will focus on $k\geq3$. From now on, \emph{we assume that $[\p]$ is a nonsingular conic}.   

\subsection{\fontfamily{cmss}\selectfont Cohomology of the twisted endomorphism bundles}\label{PushPullCalcs}

We exploit the double cover and the push-pull property of the direct image functor in order to access the cohomology of twisted endomorphism bundles of $V^\p_k$ with nonsingular $[\p]$.

\begin{proposition}\label{EndE1calc}  Assume $k\geq3$ and $d\geq0$.  If $d\geq k-1$, then\beqn\ds h^0(\CC\PP^2,\emph{End}_0V^\p_k(d)) & = & \frac{d(d+1)}{2}+(d+2)^2-k^2;\nonumber\eeqn\noin else,\beqn\ds h^0(\CC\PP^2,\emph{End}_0V^\p_k(d)) & = & \frac{d(d+1)}{2}.\nonumber\eeqn\end{proposition}

\begin{proof}Pulling back $V_k^\p$ to $\CC\PP^1\times\CC\PP^1$ gives us a surjective map\beqn(f^\p)^*V^\p_k\rightarrow\CO(0,k),\nonumber\eeqn which defines a short exact sequence\beqn0\rightarrow\CO(a,b)\rightarrow(f^\p)^*V^\p_k\rightarrow\CO(0,k)\rightarrow0.\label{SESpushpulleval}\eeqn\noin  Because of $c_1(V^\p_k)=(k-1)H$ and functoriality, we must have\beqn\CO(k-1,k-1)=\wedge^2(f^\p)^*V^\p_k=\CO(a,b+k),\nonumber\eeqn\noin and so $\CO(a,b)=\CO(k-1,-1)$. The dual sequence\beqn0\rightarrow\CO(0,-k)\rightarrow(f^\p)^*(V^\p_k)^\vee\rightarrow\CO(1-k,1)\rightarrow0\nonumber\eeqn\noin can be twisted by $\CO(d,d+k)$ to give\beqn0\rightarrow\CO(d,d)\rightarrow(f^\p)^*(V^\p_k)^\vee(d,d+k)\rightarrow\CO(d-k+1,d+k+1)\rightarrow0\nonumber\eeqn\noin  Because $H^1(\CO(d,d))=0$, we have\beqn h^0((f^\p)^*(V^\p_k)^\vee(d,d+k)) & = & h^0(\CO(d,d))+\delta_{k,d}h^0(\CO(d-k+1,d+k+1))\nonumber\\ & = & (d+1)^2+\delta_{k,d}(d+2-k)(d+2+k)\nonumber\\ & = & (d+1)^2+\delta_{k,d}((d+2)^2-k^2),\nonumber\eeqn\noin where $\delta_{k,d}=1$ if $d\geq k-1$ and $0$ otherwise.  By \eqref{SESpushpulleval}, we have\beqn h^0(\End V^\p_k(d)) & = & (d+1)^2+\delta_{k,d}((d+2)^2-k^2).\nonumber\eeqn  Removing the trace in $H^0(\CO(d))$ leaves \beqn h^0(\Endz V^\p_k(d)) & = & (d+1)^2+\delta_{k,d}((d+2)^2-k^2)-\frac{(d+1)(d+2)}{2}\nonumber\\ & = & \frac{d(d+1)}{2}+\delta_{k,d}((d+2)^2-k^2).\nonumber\eeqn\end{proof}

\noin Because stability implies $H^0(\Endz V^\p_k)=0$ for $k\geq2$, we have $H^2(\Endz V^\p_k(d))=H^0(\Endz V^\p_k(-d-3))^\vee=0$ for $d\geq0$. Combining this fact with Proposition \ref{EndE1calc} and then performing a Riemann-Roch calculation, we get:

\begin{corollary}  Again, $k\geq 3$ and $d\geq 0$.  If $d\geq k-1$, then $H^1(\CC\PP^2,\emph{End}_0V^\p_k(d))=0$; else, $h^1(\CC\PP^2,\emph{End}_0V^\p_k(d))=k^2-d^2-4d-4$.\end{corollary}

\begin{proposition}\label{PropTDim} When $k>3$, we have $h^0(\CC\PP^2,\emph{End}_0V^\p_k\tensor T)=3$.  For $k=3$, $h^0(\CC\PP^2,\emph{End}_0V^\p_3\tensor T)=8$.\end{proposition}

\begin{proof}There is another push-pull identity:\beqn H^0(\CC\PP^2,(\V)^\vee\tensor\V\tensor T) & = & H^0(\CC\PP^1\times\CC\PP^1,(f^\p)^*((\V)^\vee\tensor T)\tensor\CO(0,k)),\nonumber\eeqn\noin and so we may calculate the dimension on the right instead.  Recall from the proof of Proposition \ref{EndE1calc} the short exact sequence\beqn0\rightarrow\CO(1-k,1)\rightarrow(f^\p)^*\V\rightarrow\CO(0,k)\rightarrow0.\label{EvalSeq}\eeqn\noin  The dual sequence to (\ref{EvalSeq}) is \beqn0\rightarrow\CO(0,-k)\rightarrow (f^\p)^*(\V)^\vee\rightarrow\CO(1-k,1)\rightarrow0,\label{DualEvalSeq}\eeqn\noin from which we arrive at \beqn0\rightarrow(f^\p)^*T\rightarrow(f^\p)^*((\V)^\vee\tensor T)(0,k)\rightarrow(f^\p)^*T(1-k,k+1)\rightarrow0.\label{TwistedDualEvalSeq}\eeqn\noin  We want to calculate $H^0$ for the middle term.

At $k=2$, sequence (\ref{EvalSeq}) looks like \beqn0\rightarrow\CO(1,-1)\rightarrow(f^\p)^*V^\p_2\rightarrow\CO(0,2)\rightarrow0,\nonumber\eeqn\noin which becomes\beqn0\rightarrow\CO(2,0)\rightarrow(f^\p)^*T\rightarrow\CO(1,3)\rightarrow0\label{Twist11}\eeqn\noin after a twist by $\CO(1,1)$.  Yet another twist, this time by $\CO(1-k,1+k)$, gives\beqn0\rightarrow\CO(3-k,1+k)\rightarrow(f^\p)^*T(1-k,1+k)\rightarrow\CO(2-k,4+k)\rightarrow0.\label{Twist1-r1+r}\eeqn\noin  The cohomology of (\ref{Twist1-r1+r}) tells us that $H^0(\CC\PP^1\times\CC\PP^1,(f^\p)^*T(1-k,1+k))$ vanishes for $k>3$.  From (\ref{TwistedDualEvalSeq}), we find that\beqn H^0(\CC\PP^1\times\CC\PP^1, (f^\p)^*((\V)^\vee\tensor T)\tensor\CO(0,k)) & \cong & H^0(\CC\PP^1\times\CC\PP^1,(f^\p)^*T),\nonumber\eeqn and from (\ref{Twist11}) we can read off that $h^0(\CC\PP^1\times\CC\PP^1,(f^\p)^*T)=11$.  Traces of $T$-valued endomorphisms of $\V$ are vector fields on $\CC\PP^2$, spanning an 8-dimensional vector space.  It follows that $h^0(\CC\PP^2,\Endz\V\tensor T)=3$.

When $k=3$, sequence (\ref{Twist1-r1+r}) tells us that $H^0(\CC\PP^1\times\CC\PP^1,(f^\p)^*T(1-k,1+k))$ is not zero, but rather 5-dimensional, coming from $H^0(\CC\PP^1\times\CC\PP^1,\CO(3-k,1+k))$, and so\beqn h^0((f^\p)^*((\V)^\vee\tensor T)\tensor\CO(0,k)) & = &  h^0((f^\p)^*T)\;=\;11+5\;=\;16.\nonumber\eeqn\noin  Setting the trace to zero leaves an 8-dimensional space. \end{proof}

If we take the dual Euler sequence and twist by $\Endz\V(-3)$, we have\beqn0\rightarrow\Endz\V\tensor T^*(-3)\rightarrow(\Endz\V(-4))^{\plus 3}\rightarrow\Endz\V(-3)\rightarrow0.\nonumber\eeqn\noin  But $H^0(\Endz\V(-3))=H^2(\Endz\V)^\vee=0$ because of Proposition \ref{EndE1calc} applied to $d=0$.  It follows that $H^0(\Endz\V(-4))=0$ as well, and so in turn, we get\beqn h^2(\Endz\V\tensor T) & = & h^0(\Endz\V\tensor T^*(-3))\;=\;0.\nonumber\eeqn\noin  We use Hirzebruch-Riemann-Roch to get $h^1(\Endz\V\tensor T)=0$ when $k=3$ and $2k^2-23$ when $k>3$.

For ease of reference, the information above is summarized in the following tables.

\pagebreak

\tb{\fontfamily{cmss}\selectfont Table 1.} $k=3$\\\bcent\begin{math}\begin{array}{|c||c|c|c|}\hline \;\;\;\;\;V^\p_3\;\;\;\;\;\; & h^0 & h^1 & h^2\\\hline\hline \Endz V^\p_3 & 0 & 5 & 0\\\hline \Endz V^\p_3(1) & 1 & 0 & 0\\\hline \Endz V^\p_3(2) & 10 & 0 & 0\\\hline \Endz V^\p_3\tensor T & 8 & \;\;0\;\; & 0 \\\hline \;\;\Endz V^\p_3\tensor\wedge^2 T=\Endz V^\p_3(3)\;\; & \;\;22\;\; & 0 & \;\;0\;\;\\\hline\end{array}\end{math}\ecent

\vspace{10pt}

\tb{\fontfamily{cmss}\selectfont Table 2.} $k>3$\\\bcent\begin{math}\begin{array}{|c||c|c|c|}\hline \;\;\;\;\;V^\p_k\;\;\;\;\; & h^0 & h^1 & h^2\\\hline\hline \Endz V^\p_k & 0 & k^2-4 & 0\\\hline \Endz V^\p_k(1) & 1 & k^2-9 & 0\\\hline \Endz V^\p_k(2) & 3 & k^2-16 & 0\\\hline \Endz V^\p_k\tensor T & 3 & 2k^2-23 & 0 \\\hline \;\;\Endz V^\p_k\tensor\wedge^2 T=\Endz V^\p_k(3)\;\; & \;\;15\mbox{ if }k=4;\;6\mbox{ if }k>4\;\; & \;\;\mbox{max}(0,k^2-25)\;\; & \;\;0\;\;\\\hline\end{array}\end{math}\ecent

\vspace{10pt}

\subsection{\fontfamily{cmss}\selectfont Determining integrable Higgs fields}

We use the Euler sequence again.  Twisting the sequence by $\Endz\V\tensor$ gives the cohomology sequence\beqn0\rightarrow H^0(\Endz\V)\rightarrow\CC^3\tensor H^0(\Endz\V(1))\rightarrow H^0(\Endz\V\tensor T)\label{LESk>3}\\\rightarrow H^1(\Endz\V)\rightarrow\CC^3\tensor H^1(\Endz\V(1))\rightarrow\cdots,\nonumber\eeqn\noin in which $\CC^3$ is, as before, identified with $H^0(T(-1))$.

Exclude $k=3$ for the moment.  By stability, $H^0(\Endz\V)=0$ and so the sequence \eqref{LESk>3} begins with $\CC^3\tensor H^0(\Endz\V(1))\inject H^0(\Endz\V\tensor T)$.  On the other hand, according to Table 2, $h^0(\Endz\V(1))=1$ and $h^0(\Endz\V\tensor T)=3$, and what we are left with is\beqn H^0(\Endz\V(1))\tensor H^0(T(-1)) & \cong & H^0(\Endz\V\tensor T)\nonumber.\label{Isok>3}\eeqn\noin Since $H^0(\Endz\V(1))$ is 1-dimensional, we can fix a generator, say $\phi_0$, and then \eqref{Isok>3} says that if $\Phi\in H^0(\Endz\V\tensor T)$, then there exists a $C\in H^0(T(-1))$ such that $\Phi=\phi_0\tensor C$.   It follows that, for any $\Phi\in H^0(\Endz\V\tensor T)$, we have $\Phi\wedge\Phi=[\phi_0,\phi_0]C\wedge C=0$.  In other words, for $k>3$, every element of $H^0(\Endz\V\tensor T)$ is an integrable Higgs field.  If we choose a different nonsingular conic $[\p']$, the values in Table 2 still hold for the new bundle $V^{\p'}_k$, and so elements of $H^0(\Endz V^{\p'}_k\tensor T)$ decompose in the same manner and are integrable.  Automorphisms of any $V^\p_k$ act trivially on Higgs fields since $V^\p_k$ is stable as a bundle; that is, if $\Phi_1\neq\Phi_2$ in $H^0(\Endz\V\tensor T)$, then $(\V,\Phi_1)$ and $(\V,\Phi_2)$ are distinct points in the moduli space.  Therefore, a description of the moduli space amounts to keeping track of how $\phi$ varies with $[\p]$, and how different choices of $\phi$ and $C$ might give rise to the same $\Phi$.  In the following statement, $\Delta\subset\CC\PP^5$ denotes the set of singular conics in $\CC\PP^2$; $\mb S_k\subset\mb{CHB}\left(k-1,k(k-1)/2\right)$, the moduli space of stable co-Higgs bundles whose underlying bundles are the $k$-Schwarzenberger bundles for nonsingular $[\p]$; and $\mb S_k^*$, the complement of $\Phi=0$ in $\mb S_k$.

\begin{theorem}\label{ThmC3P5} Fix an integer $k>3$.  Then $\mb S_k^*$ is an $8$-dimensional space fibred over projective space in two ways.  The first is $\pi_A:\mb S_k^*\longrightarrow\CC\PP^2$ whose fibres $\pi_A^{-1}([C])$ are isomorphic to an unbranched double cover of $\CO_{\CC\PP^5\bsh\Delta}(-1)_0$, and the other, $\pi_B:\mb S_k^*\longrightarrow\CC\PP^5\bsh\Delta$, whose fibres are $\pi_B^{-1}([\p])\cong\CO_{\CC\PP^2}(-1)_0$; in both cases, $_0$ refers to the complement of the zero section.\end{theorem}

\begin{proof}  A Higgs field for $V^\p_k$ has the form $\Phi=\phi\tensor C$ for some $\phi\in H^0(\Endz V^\p_k(1))\cong\CC$ and some $C\in H^0(T(-1))\cong\CC^3$. If we fix the projective class of $C$, then rescaling $\phi$ changes the isomorphism class of $\Phi$, and the isomorphism classes of Higgs fields for fixed $[C]\in\CC\PP^2$ form a copy of $\CC$.  The determinant of an element $\phi$ in this line is $\alpha\p$ for some $\alpha\in\CC^*$ ($\det\phi=0$ is excluded because $\Phi\neq0$).  We also have $\det(-\phi)=\alpha\p$.  When we let $[\p]$ vary in $\CC\PP^5\bsh\Delta$, the determinants of $\phi$ form the complement of the zero section in the tautological bundle $\CO_{\CC\PP^5\Delta}(-1)$. Taken altogether, for a fixed $[C]$, the space of Higgs fields is an everywhere 2:1 cover of $\CO_{\CC\PP^5\bsh\Delta}(-1)$.   On the other hand, if we fix $[\p]$, then we get Higgs fields of the form $\phi\tensor C$ for some nonzero $\phi$ and nonzero $C$.  The Higgs field $\Phi=\phi\tensor C$ is unchanged if we replace $\phi$ with $\lambda^{-1}\phi$ and $C$ with $\lambda C$, for any $\lambda\in\CC^*$.  Under the stability condition $C\neq0$, the quotient is $\CO_{\CC\PP^2}(-1)_0$.
\end{proof}

The quotients in the preceding theorem do not depend on $k$, and so:

\begin{corollary} If $k,k'>3$, then $\mb S_k\cong\mb S_{k'}$.\end{corollary}

The analogous statement for $k=3$ is harder to ascertain because \eqref{LESk>3} reduces only so far as\beqn\;\;\;\;\;\;0\rightarrow\CC^3\tensor H^0(\Endz V^\p_3(1))\rightarrow H^0(\Endz V^\p_3\tensor T)\rightarrow H^1(\Endz V^\p_3)\rightarrow0.\label{LESk=3}\eeqn\noin The dimensions of the three terms can be read off from Table 1 as 3, 8, and 5, respectively.  The subspace $\CC^3\tensor H^0(\Endz V^\p_3(1))$ of $H^0(\Endz V^\p_3\tensor T)$ consists of integrable Higgs fields of the form $\Phi=\phi\tensor C$.  Considering only nonzero Higgs fields of this type gives us a family $\mb S_3^*$ with description identical to that in Theorem \ref{ThmC3P5}, but we cannot preclude the possibility of integrable Higgs fields that cannot be expressed as a simple tensor $\phi\tensor C$.

\section{\fontfamily{cmss}\selectfont DEFORMATIONS}

If $k>3$, then $h^1(\Endz\V)>5$ and so there are deformations of the bundle $\V$ that do not come from pushing down a rank-1 sheaf on a quadric (nonsingular or otherwise).  Is it possible to deform a point $(\V,\Phi)\in\mb S^*_k$ into a co-Higgs bundle whose underlying bundle does not come from a sheaf on a quadric? When $k=3$, $h^1(\Endz V^\p_3)=5$ and every deformation of $V^\p_3$ comes from a sheaf on a quadric, but the space of integrable Higgs fields for a fixed $V^\p_3$ is not necessarily 3-dimensional.  Can we deform a nonzero Higgs field of the form $\phi\tensor C$ for $V^\p_3$ into one that is not a simple tensor?

The answer to both questions is ``no''.  In the arguments to follow, it suffices to fix a nonzero generator $\phi_0$ of $H^0(\Endz\V(1))=\CC$.  Because $(\V,\phi_0)$ has smooth spectral cover $\CC\PP^1\times\CC\PP^1$ for $[\p]\in\CC\PP^5\bsh\Delta$, we may use in the arguments below that $\phi_0$ is regular.

\begin{theorem}\label{Thmk>=3} Fix $[\p]\in\CC\PP^5\bsh\Delta$.  For $k>3$, we have $\dim_\CC\HH^1_{(V^\p_k,\Phi)}=8$ at each $(V^\p_k,\Phi)$ with $\Phi\neq0$. For $k=3$, we have $\dim_\CC\HH^1_{(V^\p_3,\Phi)}=8$ at each $(V^\p_3,\Phi)$ for which $\Phi=\phi_0\tensor C$, where $C\in H^0(T(-1))\bsh\set{0}$.\end{theorem}

\begin{proof}  
The sequence \eqref{ESFiveTerm} contains a short exact sequence\beqn0\longrightarrow\CE^{1,0}_2\longrightarrow\HH^1_{(V,\Phi)}\longrightarrow\CE^{0,1}_2\stackrel{d_2}{\longrightarrow}\set{0}\,\subset\,\CE_2^{2,0}\nonumber,\eeqn\noin allowing us to calculate $\dim\HH^1_{(V,\Phi)}$ from $\dim\CE^{1,0}_2$ and $\dim\CE^{0,1}_2$. To see that $\mbox{im}(d_2)=\set{0}$, let $(\psi_{\alpha\beta})$ be a cocycle in\beqn\ker H^1(\Endz\V)\stackrel{[-,\Phi]}{\longrightarrow}H^1(\Endz\V\tensor T).\nonumber\eeqn\noin  Then $[\psi_{\alpha\beta},\phi_0C]=\theta_\beta C-\theta_\alpha C$, where $\theta_\alpha,\theta_\beta$ are 0-cochains for $\Endz\V\tensor T$.  But then $d_2(\psi_{\alpha\beta})\;=\;[\theta_\beta C,\phi_0C]\;=\;[\theta_\beta,\phi]\,C\wedge C\;=\;0$.\\

\tb{Case 1: $k>3$.}\\

We need to determine the dimensions of $\CE^{1,0}_2$ and $\CE^{0,1}_2$.  Directly, \beqn\CE^{1,0}_2 & = & \frac{\ker H^0(\Endz\V\tensor T)\stackrel{\wedge\Phi}{\longrightarrow}H^0(\Endz\V\tensor\wedge^2T)}{\img H^0(\Endz\V)\stackrel{\wedge\Phi}{\longrightarrow}H^0(\Endz\V\tensor T)}\nonumber\\ & = & \set{\Theta\in H^0(\Endz\V\tensor T)\;:\;\Theta\wedge\Phi=0}\nonumber\\ & = & H^0(\Endz\V\tensor T),\nonumber\eeqn\noin since $H^0(\Endz\V)=0$, and since every $\Theta\in H^0(\Endz\V\tensor T)$ can be written as $\theta\cdot C'=a\phi_0\tensor C'$ for some $a\in\CC$.  From this, we have $\Theta\wedge\Phi=\comm{\theta}{\phi_0}\cdot C'\wedge C=a\comm{\phi_0}{\phi_0}\cdot C'\wedge C$.  Therefore, $\CE_2^{1,0}=H^0(\Endz\V\tensor T)=\CC^3$.

Next, we claim that\beqn\ds\dim\CE^{0,1}_2\,:=\,\mbox{dim}\,\mbox{ker}\left(H^1(\Endz\V)\stackrel{\wedge\Phi}{\longrightarrow}H^1(\Endz\V\tensor T)\right)\,=\,5.\nonumber\eeqn\noin  This map on 1-cochains induced by $\wedge\Phi$ factors into two maps:\beqn[-,\phi_0] & : & H^1(\Endz\V)\longrightarrow H^1(\Endz\V(1)),\nonumber\eeqn\noin followed by \beqn\wedge C & : & H^1(\Endz\V(1))\longrightarrow H^1(\Endz\V\tensor T).\nonumber\eeqn  First, we show that that $\ker H^1(\Endz\V)\stackrel{[-,\phi_0]}{\longrightarrow}H^1(\Endz\V(1))$ is 5-dimensional.  

Note that since $H^0(\Endz\V(1))=\CC$, $\phi_0$ and $\phi_0^\vee$ must be scalar multiples of one another, and so the full extent of Lemma \ref{LemmRegular} applies: there exist two short exact sequences of bundles,\beqn0{\longrightarrow}\CO(-1)\stackrel{\phi_0}{\longrightarrow}\Endz\V\stackrel{[-,\phi_0]}{\longrightarrow}Q\rightarrow0\label{SESQ}\eeqn\noin and\beqn0\rightarrow Q\longrightarrow\Endz\V(1)\stackrel{\phi_0}{\longrightarrow}\CO(2)\rightarrow0\label{SESQ2}\eeqn\noin  The long exact cohomology sequence of \eqref{SESQ} has $H^0(Q)=0$ and $H^1(\Endz\V)\cong H^1(Q)$.  Using these facts, \eqref{SESQ2} has cohomology\beqn0\rightarrow H^0(\Endz\V(1))\rightarrow H^0(\CO(2))\rightarrow H^1(\Endz\V)\stackrel{[-,\phi_0]}{\longrightarrow}H^1(\Endz\V(1))\rightarrow0,\nonumber\eeqn\noin which gives us the required surjectivity.  We can also read off from the sequence that $\ker[-,\phi_0]:H^1(\Endz\V)\rightarrow H^1(\Endz\V(1))$ is 5-dimensional.

Now we show that the second map, $H^1(\Endz\V(1))\stackrel{\wedge C}{\longrightarrow}H^1(\Endz\V\tensor T)$, is injective.  Note that the exact sequence\beqn0\rightarrow\Endz\V(1)\rightarrow\Endz\V\tensor T\rightarrow\Endz\V(2)\tensor\MF I_x\rightarrow 0\nonumber\eeqn\noin coming from the map $\CO\rightarrow T(-1)$ given by $f\mapsto fC$.  Since $C$ is not identically zero, it vanishes only at a point $x\in\CC\PP^2$, which defines $\MF I_x$.  The long cohomology sequence begins with the left-exact sequence\beqn0\rightarrow H^0(\Endz\V(1))\rightarrow(\Endz\V\tensor T)\rightarrow H^0(\Endz\V(2)\tensor\MF I_x)\nonumber\eeqn\noin in which the first term is 1-dimensional and the second is 3-dimensional.  Referring to the cohomology table we also know that $h^0(\Endz\V(2))=3$, and so the constraint that sections vanish at $x$ means $h^0(\Endz\V(2)\tensor\MF I_x)=2$.  This makes the left-exact sequence fully exact.  Therefore, $H^1(\Endz\V(1))\rightarrow H^1(\Endz\V\tensor T)$ is injective.

Hence $\dim\CE_2^{0,1}=5$, and so $\dim\HH^1_{(\V,\Phi)}=8$.\\

\tb{Case 2: $k=3$.}\\

We claim that $\dim\CE^{1,0}_2=3$.  Again, we have\beqn\CE_2^{1,0} & = & \set{\Theta\in H^0(\Endz V^\p_3\tensor T)\;:\;\Theta\wedge\Phi=0}.\nonumber\eeqn\noin  Recall also the short exact sequence\beqn0\rightarrow\CC^3\tensor H^0(\Endz V^\p_3(1))\rightarrow H^0(\Endz V^\p_3\tensor T)\rightarrow H^1(\Endz V^\p_3)\rightarrow0.\nonumber\eeqn\noin We know $\CC^3\tensor H^0(\Endz V^\p_3(1))\subset\CE_2^{1,0}$.  We want to show that the inclusion is an inequality.

Use $x$ again for the vanishing point of $C$; $\MF I_x$, the ideal sheaf concentrated there.  The equation $\Psi\wedge\Phi=0$ can be written $[\phi_0,\Theta\wedge C]=0$.  Solving this and extending via Hartogs over $x$, we get $\Theta=s\phi_0$, where $s\in H^0(\CO(1))$ because $\Theta\in H^0(\Endz V^\p_3(2)\tensor\MF I_x)$. Since $\phi_0$ is regular, it must be $s$ that vanishes at $x$; in particular, $s$ passes through $x$.  This means that there are two degrees of freedom in choosing $s$: a single restriction applied to $H^0(\CO(1))=\CC^3$.  Consider now the map on functions given by $f\mapsto Cf$.  This gives rise to an exact sequence of sheaves\beqn0\rightarrow\CO\rightarrow T(-1)\rightarrow\MF I_x\tensor\CO(1)\rightarrow0,\nonumber\eeqn\noin which in turn gives us\beqn0\rightarrow\Endz V^\p_3(1)\rightarrow\Endz V^\p_3\tensor T\stackrel{\wedge C}{\longrightarrow}\Endz V^\p_3(2)\tensor\MF I_x\rightarrow0\nonumber\eeqn\noin once we apply $\Endz V^\p_3(1)\tensor$.  Applying $H^0$ and noting $H^1(\Endz V^p_3(1))=0$, we have\beqn\;\;\;\;\;\;\;\;\;0\rightarrow H^0(\Endz V^\p_3(1))\rightarrow H^0(\Endz V^\p_3\tensor T)\stackrel{\wedge C}{\longrightarrow}H^0(\Endz V^\p_3(2)\tensor\MF I_x)\rightarrow0,\nonumber\eeqn\noin in which the first space, $H^0(\Endz V^\p_3(1))$, is 1-dimensional.  The problem is now about determining which elements of $H^0(\Endz V^\p_3\tensor T)$ go to elements of the form $s\phi_0$ in $H^0(\Endz V^\p_3(2)\tensor\MF I_x)$.  Since such elements form a 2-dimensional subspace of $H^0(\Endz V^\p_3(2)\tensor\MF I_x)$, and since the kernel of the exact sequence is 1-dimensional, we conclude that inside $H^0(\Endz E_3\tensor T)$ is a 3-dimensional subspace whose elements take the desired form after $\wedge C$, and this subspace is precisely $\CC^3\tensor H^0(\Endz V^\p_3(1))$.

To finish, we need to show that\beqn\dim\CE^{0,1}_2\;:=\;\ker H^1(\Endz V^\p_3)\stackrel{\wedge\Phi}{\longrightarrow}H^1(\Endz V^\p_3\tensor T)\;=\;5.\nonumber\eeqn\noin This follows straight away from \beqn h^1(\Endz V^\p_3)=5\;\mbox{ and }\;h^1(\Endz V^\p_3\tensor T)=0,\nonumber\eeqn\noin as listed in the cohomology table.  Therefore, $\dim\HH^1_{(V^\p_3,\Phi)}=8$.

\end{proof}

\begin{remark}
Calculations along the lines of those in the proof of Theorem \ref{Thmk>=3} show that, for $k\geq3$, $\dim\HH^2\neq0$ at any $(V^\p_k,\phi\tensor C)$ with $\phi\tensor C\neq0$.  Although the obstruction space is large (it varies with $k$, but is at least $9$-dimensional), it is clear from Theorem \ref{ThmC3P5} that, exluding zero Higgs fields, the co-Higgs bundles of our construction are situated in a smooth subvariety of their respective moduli spaces.  This subvariety is isomorphic to $\mb S_k^*$.
\end{remark}

Theorem \ref{Thmk>=3}, when combined with earlier results for $\CO\plus\CO(-1)$, $\CO\plus\CO$, and $T$, leads us to the following set of conclusions:

\begin{corollary}\label{CorRigid} Let $k$ be a nonnegative integer.\bitem\item[\tb{(i)}] If $k=0,1,2$, let $[\p]\in\CC\PP^5$ and $\Phi$ be any stable Higgs field for $\V$ satisfying $\Phi\wedge\Phi=0$. 
\item[\tb{(ii)}] If $k=3$, let $[\p]\in\CC\PP^5\bsh\Delta$ and let $\Phi$ be any Higgs field for $V^\p_3$ of the form $\phi\tensor C$ for some $\phi\in H^0(\CC\PP^2,\emph{\mbox{End}}_0V^\p_3(1))\bsh\set{0}$ and some $C\in H^0(T(-1))\bsh\set{0}$. \item[\tb{(iii)}] If $k>3$, let $[\p]\in\CC\PP^5\bsh\Delta$ and let $\Phi$ be any nonzero Higgs field for $\V$.\eitem\noin  Then $(\V,\Phi)$ can only be deformed to a co-Higgs bundle whose underlying bundle is $V^{\p'}_k$, where $[\p']$ is a (possibly singular) conic in $\CC\PP^2$. If a first-order deformation $(V',\Phi')$ of $(V^\p_3,\phi\tensor C)$ has underlying bundle $V'=V^{\p'}_3$ for some $[\p']\in\CC\PP^5\bsh\Delta$, then $\Phi'=\phi'\tensor C'$ for some $\phi'\in H^0(\CC\PP^2,\emph{\mbox{End}}_0V^{\p'}_3(1))$ and some $C'\in H^0(T(-1))$.\end{corollary}

\begin{remark}
Corollary \ref{CorRigid} says that, nearby to a co-Higgs bundle of Schwarzenberger type, there are no co-Higgs bundles of a different type.  For $k\geq3$, corollary \ref{CorRigid} does not necessarily imply that the families $\mb S_k^*$ are topological components of their respective moduli spaces of co-Higgs bundles.  Because our arguments rely on regularity, we have not ruled out the following possibilities: deforming out of the Schwarzenberger family from a point $(\V,\Phi)$ for which $[\p]$ is singular; or, since the zero Higgs field is stable in these cases, constructing a connected path from a point $(\V,\Phi)$ to the zero Higgs field and then into the moduli space of stable bundles on $\CC\PP^2$, which is a locus in the moduli space of stable co-Higgs bundles.
\end{remark}

\section{\fontfamily{cmss}\selectfont REMARKS}\label{Further}

\subsection{\fontfamily{cmss}\selectfont Spectral interpretation and Hitchin map}

In our construction, we fashion co-Higgs bundles from direct images of line bundles on $\CC\PP^1\times\CC\PP^1$.  But $\CC\PP^1\times\CC\PP^1$ is a naturally a subvariety of $\CC\PP^3$, via the Segre embedding, and $\CC\PP^3$ is the one-point compactification of the total space of $\CO(1)\rightarrow\CC\PP^2$.  By the spectral correspondence for Higgs bundles, this means that pushing down line bundles from the quadric yields not only vector bundles on $\CC\PP^2$, but also $\CO(1)$-valued Higgs fields, which we have already seen appearing as intermediate objects in our construction.

We can view the entire $\CO(1)$-valued side of the construction as occurring within $\CC\PP^3$.  Let $u,v,w$ be affine coordinates on $\CC\PP^3$.  We fix the plane $P=\CC\PP^2$ defined by $w=0$ and a smooth quadric $Q_\p=\CC\PP^1\times\CC\PP^1$ defined by $w^2-\p(u,v)=0$, where $\p\in H^0(P,\CO(2))$.  We can project $Q_\p$ onto $P$ from the point at infinity, giving us a projection map $f^\p$ whose branch locus is the conic determined by $\p$.  The line bundle $\CO(1,1)\rightarrow Q_\p$, which is the pullback of $\CO(1)$ from $\CC\PP^3$, has a four-dimensional space of sections.  Three independent generators for this space are generators of $H^0(P,\CO(1))\cong\CC^3$, pulled back to $Q_\p$.  The fourth generator is the tautological section $s$ of the pullback of $\CO(1)\rightarrow P$ to its own total space.  This section gives us a multiplication $\CO(0,k)\stackrel{s}{\rightarrow}\CO(1,k+1)$ over $Q_\p$ that can be pushed down via the projection map $f^\p$, giving a twisted endomorphism $\phi=f^\p_*s:V^\p_k\rightarrow V^\p_k(1)$ with determinant $\p\in H^0(P,\CO(2))$. 

Whenever we have $\phi\in H^0(\Endz\V(1))$ with determinant $\p$, we can get a $T$-valued Higgs field for $\V$ by tensoring $\phi$ with a section $C$ of $T(-1)$.  The determinant transforms, accordingly, from $\p$ to $\p C\tensor C$.  The spectral variety $Q_\p$ transforms into a subvariety $Q_{\p,C}\subset\mbox{Tot}(\pi:T\rightarrow\CC\PP^2)$ cut out by the equation $\eta^2-\pi^*(\p C^{\tensor2})=0$, where $\eta$ is the tautological section of $\pi^*T$. Since $C$ vanishes at a point in $\CC\PP^2$, $Q_{\p,C}$ will be a singular subvariety.  This is consistent with the fact that every Higgs field we have constructed vanishes at a point in $\CC\PP^2$, because of the factor of $C$, and therefore $\Phi$ could not been pushed down from a smooth spectral cover in $\mbox{Tot}(T)$.

%Finally, note that because the fiber of the Hitchin map in the $O(1)$-valued side of the construction is discrete (it is a Jacobian of $\CC\PP^1\times\CC\PP^1$), so too is the fiber on the $T$-valued side.

\subsection{\fontfamily{cmss}\selectfont Extending across $\Delta$}

When $\p$ is reducible, the associated double cover of $\CC\PP^2$ is a singular quadric. The line bundle $\CO(0,k)$ is replaced by a reflexive sheaf $I$.  The direct image of $I$ under $f^\p$ is a reflexive sheaf on $\CC\PP^2$, and therefore is a rank-2 vector bundle. For $k\geq3$, a natural question is whether theorems \ref{ThmC3P5} and \ref{Thmk>=3} can be extended across $\Delta$.  We leave this question for future work, but offer some basic comments.  The functions $h^{i,k}(\p):=\dim H^i(\Endz\V\tensor\wedge^kT)$ are constant over the nonsingular locus $\CC\PP^5\bsh\Delta$.  We do not have a unique extension theorem for the $h^{i,k}$.  Alternatively, we could try to calculate the numbers $h^{i,k}$ over $\Delta$ directly.  The push-pull functoriality does not depend on the smoothness of the double cover.  In particular, since $f^\p$ is always finite, there is still an exact sequence of functors\beqn0\rightarrow\CO(-R)\tensor j^*\rightarrow (f^\p)^* f^\p_*\rightarrow\mbox{id}\rightarrow0,\nonumber\eeqn\noin where $R$ is the ramification divisor in the quadric and $j$ is the sheet interchange.  However, the arguments used to prove Theorem \ref{Thmk>=3} rely on the regularity of $\phi$, which is unavailable in the singular case.

%%%
%%%
%%%
%%%
%%%%%%%%%%%%%%%%%%%%%%%%%%%%%%%%%%%%%%%%%%%%%%%%%%%%%%%%%%%%%%%%%%%%%%%%%%%%%%%%%%%%%%%%%%%%%%%%%%%%
%%%%%%%%%%%%%%%%%%%%%%%%%%%%%%%%%%%%%%%%%%%%%%%%%%%%%%%%%%%%%%%%%%%%%%%%%%%%%%%%%%%%%%%%%%%%%%%%%%%%
%%%%%%%%%%%%%%%%%%%%%%%%%%%%%%%%%%%%%%%%%%%%%%%%%%%%%%%%%%%%%%%%%%%%%%%%%%%%%%%%%%%%%%%%%%%%%%%%%%%%

\bibliographystyle{acm}	% (uses file "acm.bst")
\bibliography{biblio-StevenRayan-24Sept2013}

\begin{thebibliography}{10}

\bibitem{BR:94}
{\sc Biswas, I., and Ramanan, S.}
\newblock An infinitesimal study of the moduli of {H}itchin pairs.
\newblock {\em J. London Math. Soc. (2) 49}, 2 (1994), 219--231.

\bibitem{BS:92}
{\sc Bohnhorst, G., and Spindler, H.}
\newblock The stability of certain vector bundles on {${\bf P}\sp n$}.
\newblock In {\em Complex {A}lgebraic {V}arieties ({B}ayreuth, 1990)},
  vol.~1507 of {\em Lecture Notes in Math.} Springer, Berlin, 1992, pp.~39--50.

\bibitem{FB:95}
{\sc Bottacin, F.}
\newblock Symplectic geometry on moduli spaces of stable pairs.
\newblock {\em Ann. Sci. \'Ecole Norm. Sup. (4) 28}, 4 (1995), 391--433.

\bibitem{SB:13}
{\sc Boucksom, S., Demailly, J.-P., P{\u{a}}un, M., and Peternell, T.}
\newblock The pseudo-effective cone of a compact {K}\"ahler manifold and
  varieties of negative {K}odaira dimension.
\newblock {\em J. Algebraic Geom. 22}, 2 (2013), 201--248.

\bibitem{RD:95}
{\sc Donagi, R.}
\newblock Spectral covers.
\newblock In {\em Current {T}opics in {C}omplex {A}lgebraic {G}eometry
  ({B}erkeley, {CA}, 1992/93)}, vol.~28 of {\em Math. Sci. Res. Inst. Publ.}
  Cambridge Univ. Press, Cambridge, 1995, pp.~65--86.

\bibitem{RF:98}
{\sc Friedman, R.}
\newblock {\em Algebraic {S}urfaces and {H}olomorphic {V}ector {B}undles}.
\newblock Springer Universitext, New York, N.Y., 1998.

\bibitem{MG:07}
{\sc Gualtieri, M.}
\newblock Branes on {P}oisson varieties.
\newblock In {\em The {M}any {F}acets of {G}eometry: {A} {T}ribute to {N}igel
  {H}itchin}. OUP, Oxford, 2010, pp.~368--394.

\bibitem{MG:11}
{\sc Gualtieri, M.}
\newblock Generalized complex geometry.
\newblock {\em Ann. of Math. (2) 174}, 1 (2011), 75--123.

\bibitem{NJH:86}
{\sc Hitchin, N.~J.}
\newblock The self-duality equations on a {R}iemann surface.
\newblock {\em Proc. London Math. Soc. (3) 55}, 1 (1987), 59--126.

\bibitem{NJH:10IIa}
{\sc Hitchin, N.~J.}
\newblock Generalized holomorphic bundles and the {$B$}-field action.
\newblock {\em J. Geom. Phys. 61}, 1 (2011), 352--362.

\bibitem{KP:94}
{\sc Katzarkov, L., and Pantev, T.}
\newblock Representations of fundamental groups whose {H}iggs bundles are
  pullbacks.
\newblock {\em J. Differential Geom. 39}, 1 (1994), 103--121.

\bibitem{KW:70}
{\sc Kobayashi, S., and Wu, H.-H.}
\newblock On holomorphic sections of certain hermitian vector bundles.
\newblock {\em Math. Ann. 189\/} (1970), 1--4.

\bibitem{NN:91}
{\sc Nitsure, N.}
\newblock Moduli space of semistable pairs on a curve.
\newblock {\em Proc. London Math. Soc. (3) 62}, 2 (1991), 275--300.

\bibitem{TP:11}
{\sc Peternell, T.}
\newblock Generically nef vector bundles and geometric applications.
\newblock In {\em Complex and {D}ifferential {G}eometry: {C}onference {H}eld at
  {L}eibniz {U}niversit{\"a}t Hannover, September 14--18, 2009\/} (Berlin,
  2011), Springer-Verlag, pp.~175--189.

\bibitem{SSR:10}
{\sc Rayan, S.}
\newblock Co-{H}iggs bundles on {$\PP^1$}.
\newblock ar{X}iv:math/1010.2656 [math.{AG}], 2010.

\bibitem{RT:11}
{\sc Ross, J., and Thomas, R.}
\newblock Weighted projective embeddings, stability of orbifolds, and constant
  scalar curvature {K}\"ahler metrics.
\newblock {\em J. Differential Geom. 88}, 1 (2011), 109--159.

\bibitem{RLES:61I}
{\sc Schwarzenberger, R. L.~E.}
\newblock Vector bundles on the projective plane.
\newblock {\em Proc. London Math. Soc. (3) 11\/} (1961), 623--640.

\bibitem{CSS:67}
{\sc Seshadri, C.~S.}
\newblock Space of unitary vector bundles on a compact {R}iemann surface.
\newblock {\em Ann. of Math. (2) 85\/} (1967), 303--336.

\bibitem{TSCH:94I}
{\sc Simpson, C.~T.}
\newblock Moduli of representations of the fundamental group of a smooth
  projective variety; {I}.
\newblock {\em Inst. Hautes \'Etudes Sci. Publ. Math.}, 79 (1994), 47--129.

\bibitem{TSCH:94II}
{\sc Simpson, C.~T.}
\newblock Moduli of representations of the fundamental group of a smooth
  projective variety; {II}.
\newblock {\em Inst. Hautes \'Etudes Sci. Publ. Math.}, 80 (1994), 5--79.

\end{thebibliography}

\end{document}